\documentclass[11pt]{amsart}

\usepackage{etoolbox}

\makeatletter
\patchcmd{\@setauthors}{\MakeUppercase}{}{}{}
\patchcmd{\section}{\scshape}{}{}{}
\makeatother

\usepackage[utf8]{inputenc}
\usepackage[T1]{fontenc}
\usepackage[english]{babel}
\usepackage{fullpage}

\usepackage{amssymb,libertine}

\usepackage{array}
\usepackage{bbm}

\usepackage{mathtools}
\usepackage[svgnames]{xcolor}
\usepackage{hyperref}
\usepackage{comment}
\usepackage[colorinlistoftodos]{todonotes}

\usepackage{enumitem}
\usepackage{wrapfig}
\usepackage{tikz}
\usepackage{color}
\usepackage{graphicx}
\usepackage{amssymb}
\usepackage{epstopdf}
 \usepackage{latexsym}
\usepackage{amsmath}
\usepackage{amsfonts} 
\usepackage{amsthm}
\usepackage[export]{adjustbox}
\def\lap{\mathcal{L}}
   \def\CC{\mathbb{C}}
    \def\DD{\mathbb{D}}
    \def\NN{\mathbb{N}}
    
    \def\RR{\mathbb{R}}
    \def\ZZ{\mathbb{Z}}
    \newtheorem{Proposition}{Proposition}
\newtheorem{Theorem}[Proposition]{Theorem}
\newtheorem{Lemma}[Proposition]{Lemma}
\newtheorem{Definition}[Proposition]{Definition}
\newtheorem{Corollary}[Proposition]{Corollary}
\newtheorem{Remark}[Proposition]{Remark}
\newtheorem{Note}[Proposition]{Note}
  
\def\be{\begin{equation}}
\def\ee{\end{equation}}
    
\def\ge{\geqslant}
\def\le{\leqslant}

\def\bd{\begin{Definition}}
\def\ed{\end{Definition}}
\def\bt{\begin{Theorem}}
\def\et{\end{Theorem}}

\def\epsilon{\varepsilon}
\def\bel{\begin{equation}\label}
\def\ee{\end{equation}}

\def\Ei\text{Ei}

\usepackage{mathrsfs}
\def\phi{\varphi}
\title{Global Rational approximations of functions with factorially divergent asymptotic series}
\author{N. Castillo, O. Costin and  R.D. Costin}

\begin{document}
\today 

\maketitle

\begin{abstract}

We construct a new type of convergent, and  asymptotic, representations, dyadic expansions. Their convergence is geometric and the region of convergence often extends from infinity down to $0^+$.  We show that dyadic expansions are numerically efficient representations.

For
special functions such as Bessel, Airy, Ei, erfc, Gamma, etc. the region of convergence of dyadic series is the complex plane minus a ray, with this cut chosen at will. Dyadic expansions thus provide uniform, geometrically convergent asymptotic expansions including near antistokes rays. 

We prove that relatively  general functions,  \'Ecalle resurgent ones, possess convergent  dyadic expansions. 

These expansions extend to operators, resulting in representations of the resolvent of self-adjoint operators as series in  terms of the associated unitary evolution operator evaluated at some prescribed discrete times  (alternatively, for positive operators, in terms of the generated semigroup).

\end{abstract}

\section{Introduction} 

\subsection{Classical approximations.}
\subsubsection{Functions given by a convergent power series}
Rational approximations are a powerful tool for generating efficient approximations for functions specified by a given convergent series, often beyond the radius of convergence of this series.  One notable example are the Pad\'e approximants, which converge geometrically (in capacity).  
The domain of convergence of a Pad\'e expansion is dictated by the structure and nature of the singularities of the function $f$ given by the series.  For instance, if $f$ has only one singularity, a branch point, then Pad\'e approximants will place their poles along a ray emanating from the branch point and going to infinity. There are libraries of classical references for the theory of Pad\'e expansions such as \cite{Wall}.  The fundamental paper by Stahl \cite{Stahl} contains powerful and detailed results about Pad\'e convergence.  

\subsubsection{Functions with divergent power series}
In many applications however, equations can only be solved by divergent series. These formal expansions are often known, from general theory, to be asymptotic to actual solutions. Let us place the asymptotic limit \emph{conventionally at $+\infty$}. When these asymptotic series 
 are Borel summable in a strip containing the non-negative real axis, there exist rising factorial expansions (factorial series, Horn series) converging in a half-plane to the Borel sum. The Borel sum is often guaranteed (by theorems) to be an actual solution of the problem of origin.

A classical rising factorial expansion for large $x$ in the open right half-plane is a series of the form 
$\displaystyle {S=  \sum_{k=1}^{\infty}\frac{c_k}{(x)_k}  }$
 where
  \begin{equation}\label{Poch}
 (x)_k:=x(x+1)\cdots (x+k-1)=\frac{\Gamma(x+k)}{\Gamma(x)}
  \end{equation}
 is known as the Pochhammer symbol, or rising factorial. 

Factorial series have a long history going back to {Stirling}, {Jensen}, {Landau}, {N\"orlund} and {Horn} (see, e.g. \cite{Stirling}, \cite{Jensen}, \cite{Landau}, \cite{N\"orlund}, \cite{Horn}).  Excellent {introductions} to the classical theory of factorial series and their application to solving ODEs  can be found in the books by N\"orlund \cite{N\"orlund} and Wasow \cite{Wasow}; see also \cite{Paris} Ch.4.

Note that since  $(x)_{k+1}$ behaves roughly like $k!$ for large $k$, then if the $c_k$ grow at most like $k!$ then the series $S$ converges even when its {\em asymptotic series} in powers of $1/x$ has empty domain of convergence; we elaborate more on this phenomenon in \S\ref{Classical}.  

Recent use of factorial expansions to tackle  divergent perturbation series in quantum mechanics and quantum field theory (see e.g. \cite{Jen}) triggered considerable renewed interest and substantial  literature. An excellent account of new developments is \cite{Weniger2000}; see  also \cite{Dunster,Adri,We97,JW} and references therein. Factorial series also play a major role in the use of sequence transformations in optics, see \cite{We89} \cite{We97}.

\subsection{Limitations of classical factorial expansions} \label{Limitations}
Most often, the classical factorial expansions used in ODEs and physics have two major limitations: (1) slow convergence, at best power-like, and (2) a limited domain of convergence (usually unrelated to the functions represented): a half-plane. The boundary of this half-plane is separated by a positive angular distance from the important antistokes rays (where the transition between power-like decay and oscillatory behavior occurs); this angular separation is necessary: see \cite{Wasow}, Theorem 46.2, p. 329 combined with the fact that, with the normalization in Wasow, Borel summability fails along $\RR^+$, see \cite{Duke} Theorem 1 (i).

\subsection{Overview of the paper}  
Throughout this work we shall use the term {\em{dyadic series}} to refer to series associated to binary (i.e. $2^k$) partitions, such as the double sum \eqref{dyadicf}.  The geometric nature of the gaps ensures geometric convergence of the ensuing series. 
For functions whose inverse Laplace transform has only one singularity on the first Riemann sheet\footnote{This is the case for many   
classical special functions \cite{Book}.}, 
the singularities of the dyadic expansions accumulate along a ray $R$.  In this case, the domain of geometric convergence of our expansions is $D=\CC\setminus R$. $R$ can be placed arbitrarily in a closed quadrant between the Stokes line and an antistokes line; this allows for providing expansions convergent in the full sector where solutions have asymptotic series, as well as in the region with the classical Stokes phenomenon, capturing the transition between an asymptotic series with a monotonic behavior and another one with oscillatory behavior 
(see examples in \S\ref{Eidyadic}). This is not possible with classical factorial series.

As in the case of Pad\'e approximants or factorial series, dyadic expansions can be calculated based upon the asymptotic series alone.

We then extend our theory to operators and develop dyadic resolvent decompositions for self-adjoint operators in terms of the associated unitary evolution, and, for positive operators, in terms of the evolution semigroup.

Finally, we address the problem of representing functions with several singularities by developing  a general theory of decomposition  of functions into  simpler ones, {\em{ function} elements}
, defined in \S\ref{FuncEl}.

The content is as follows. We start with an overview of classical factorial series and their convergence, which will motivate and clarify our approach, \S\ref{Classical}.

Section \ref{S3} contains the main results and techniques, which are subsequently generalized and applied. The  dyadic series \eqref{dyadicf} for functions $f$ which are Laplace transforms of resurgent elements are proved to converge in a cut plane. Estimates of the remainder in various regimes are given in Theorem\,\ref{ElementDecomp}.

In Section \ref{Eidyadic} we obtain dyadic  expansions of various special functions. While Theorem\,\ref{ElementDecomp} can be applied, sometimes it is easier to obtain the expansion directly, and we explain how. The fist two examples  concern the exponential integral Ei.  The next examples are the dyadic expansions for Airy functions in \S\ref{Airy}, and  for general Bessel functions in \S\ref{genBessel}.
In \S\ref{Section5} we address the question of calculating, in a practical way,  the coefficients of dyadic expansions for more general functions. 
In \S\ref{Psi} we develop the dyadic series for the Psi function. As a consequence we obtain the identity \eqref{eq:ide}, which seems to be new.
In \S\ref{dupl} we find interesting connections between dyadic expansions of the Lerch Phi function and polylogs. We use these connections to obtain interesting function identities \eqref{EiLer}--\eqref{eq:resol2}.
In \S\ref{op} we find dyadic resolvent representations of self-adjoint operators in terms of the unitary evolution operator at specific discrete times, \eqref{eq:resol2}.

In \S\ref{Resfun}  we address the general case: we develop the theory of constructing geometrically convergent dyadic expansions  for typical \'Ecalle resurgent functions. Since, by definition,  resurgent divergent series are \'Ecalle-Borel summable (to  resurgent functions, cf. footnote 1),   such  series  are also  resummable in terms of dyadic expansions.
Our theory extends naturally to transseriable functions, but we do not pursue this in the present paper.

The dyadic series introduced here are new types of representations. As it is often the case, new representations can be used to obtain new type of identities. For example, we obtain here: 
\begin{itemize}
\item the exponential integral can be written as a convergent series of Lerch $\Phi$ functions:
\begin{equation}
    \label{EiLer}
     \mathrm{e}^{-x}{\rm Ei}^+(x)=-\Phi(-1,1,\tfrac x{i\pi})+\sum_{k=1}^\infty\Phi(- \mathrm{e}^{\pi i/2^k},1,2^k \tfrac x{i\pi})\ \ \ \text{ for }x\in\CC\setminus i(-\infty,0]
\end{equation}
\item also
\begin{equation}
  \label{eq:EiPhi}
\mathrm{e}^y{\rm Ei}^+(-y)=\Phi(\mathrm{e}^{-1},1,y)-\sum_{k=1}^\infty\Phi(- \mathrm{e}^{-2^{-k}},1,2^k y),\ \ \ \text{ for }|\arg\,y|<\pi/2
\end{equation}
\item the $\Psi$ function satisfies the identity
\begin{equation}
  \label{eq:ide}
  \Psi \left( x+1 \right) =\ln  x  -\frac12 \sum _{k=0}^{
\infty }\left[\Psi \left( {2}^{k}x+1 \right) -\Psi \left( {2}^{k}x+\tfrac12\right)\right]
\end{equation}

\item representations of the resolvent of a self-adjoint operator in a series involving the unitary evolution operator at specific discrete times:
 \begin{equation}
  \label{eq:resol2}
(A-i\lambda)^{-1}  =i\sum_{j=0}^{\infty}\mathrm{e}^{-j\lambda}U_{j} -i\lim_{\ell\to\infty}\sum_{k=1}^{\ell}\sum_{j=0}^{\infty}(-1)^j2^{-k}\mathrm{e}^{-{j\lambda}/{2^k}}U_{j2^{-k}}
\end{equation}

\end{itemize}

\subsection{Why do factorial expansions usually converge in only a half-plane?  Why is convergence only power-like? }\label{Classical}
In this section we look at typical factorial expansions. We contrast them with a special case, the factorial expansions for the Lerch $\Phi$ transcendent, which converges geometrically on the full domain of analyticity. This contrast will clarify our method of eliminating the limitations noted in \S\ref{Limitations}.

The connection of factorial expansions to Borel summation was made already in \cite{N\"orlund}. Assume $f$ is the Borel sum of a series, that is, $f$ is the Laplace Transform of a function $F$
\begin{equation}\label{fisLF}
f(x)=\int_0^\infty F(p)\mathrm{e}^{-px}\, dp:=(\mathcal{L}F)(x)
\end{equation}
 where $F$ is analytic in an open sector containing $\RR^+$, analytic at $p=0$ and exponentially bounded at infinity. The asymptotic series of $f$ for large $\Re x$ is related to the Mclaurin series of $F$: this follows from Watson's lemma \cite{Wasow} or, in this case, simply by integration by parts: for $x$ large enough we have
\begin{equation}
  \label{serfF}
  f(x)=\frac1xF(0)+\frac1{x^2}F'(0)+\cdots+\frac1{x^{n}}F^{(n-1)}(0)+\frac1{x^{n}}\int_0^\infty F^{(n)}(p)\mathrm{e}^{-px}\, dp
\end{equation} 
Integration by parts results in a growing power of $\frac{d}{dp}$ and thus, by Cauchy's theorem, leads to factorial divergence of the asymptotic series of $f$, unless $F$ is entire (rarely the case in applications). 

N\"orlund notices however that the simple change of variables $\phi(s)=F(-\ln s)$ brings the representation \eqref{fisLF} of $f$  to the form 

\begin{equation}
  \label{eq:Mellin}
  f(x)=\int_0^1s^{x-1}\phi(s)ds
\end{equation}
Now integration by parts gives the {\em factorial} expansion 

\begin{equation}
  \label{eq:fct2}
 f(x)= \phi(1)\, \frac{1}x-\phi'(1)\, \frac{1}{(x)_2}+\cdots +\frac{(-1)^{n-1}}{(x)_n}\,\phi^{(n-1)}(1) +\frac{(-1)^n}{(x)_n}\int_0^1 s^{x+n-1}\phi^{(n)}(s)ds
\end{equation}
or, without remainder, we have the factorial series, (a formal series, for now)
\begin{equation}
  \label{eq:fct2sinf}
\tilde{f}(x)= {\sum_{k=0}^{\infty}(-1)^k\frac{\phi^{(k)}(1)}{(x)_{k+1}}}
\end{equation}
\begin{Note}{\rm 
Since $F$ is analytic at zero, $\phi$ is analytic at one. Using Stirling's formula in \eqref{Poch}, we see that, for large $k$,  the $(k+1)$'th term of the expansion \eqref{eq:fct2sinf} behaves like 
\begin{equation}\label{ktermphi}
  (-1)^k\, {\Gamma(x)} \,\frac{ \phi^{(k)}(1)\ }{k!}\, k^{-x}  
\end{equation}
Due to the $1/k!$ factor in \eqref{ktermphi} the series $\tilde{f}(x)$ can converge even if the power series obtained from \eqref{serfF} is factorially divergent.  
 }
\end{Note}

\begin{Note}\label{Note2}
  {\rm
    For $\tilde{f}$  to converge, {\eqref{ktermphi} shows that} $\phi$ needs to be analytic in ${\DD_{1}(1)}$, the disk of radius one centered at {$s=1$}. Indeed, $\sum_k \frac{ \phi^{(k)}(1)\ }{k!}(-1)^k$ is the Taylor series of $\phi$ about $1$, evaluated at $s=0$.  Furthermore if $\tilde{f}$ converges  geometrically, then this implies analyticity of $\phi$ in a disk larger than $\overline{\DD_{1}(1)}$.
 
  In applications $\phi$ is often singular at $s=0$. Such a singularity is one source of the limitations of factorial expansions. To illustrate this on a concrete example, we express the exponential integral $ {\rm Ei} (x):={\rm{e}}^{-x}\mathcal{L}(\tfrac1{1+p})(x)$ 
  in the form  \eqref{eq:Mellin} and obtain
\begin{equation}
    \label{eq:Ei01}
     {\rm Ei} (x)={\rm e}^{-x}\int_0^1 \frac{s^{x-1}}{1-\ln(s)}ds=:{\rm e}^{-x}\int_0^1 s^{x-1}\phi(s)ds
  \end{equation}
  The presence of the logarithmic singularity in \eqref{eq:Ei01} shows that analyticity in $\overline{\DD_{1}(1)}$ is not satisfied.

  }\end{Note}
We next examine the connection between half-plane convergence of factorial series and Borel summability. Note first that, if $\mathcal{L}F=f$, then  $\mathcal{L}[(1-e^{-p})F(p)](x)=f(x)-f(x+1)$, from which it follows immediately that, if $k+1\in\NN$, then
   \begin{equation}
     \label{eq:identexp}
     [\mathcal{L}(1-e^{-p})^k](x)=\frac{k!}{(x)_{k+1}}
   \end{equation}
    Hence, $F$, the formal inverse Laplace transform of $\tilde{f}$ in \eqref{eq:fct2sinf} is the function series
   \begin{equation}
     \label{eq:expseries}
     \sum_{k=0}^\infty \frac{(-1)^k\phi^{(k)}(1)}{k!}(1-e^{-p})^k
   \end{equation}

\ 

We now contrast the slow convergence of typical factorial expansions with the convergence of the factorial expansions of the Lerch $\Phi$ transcendent, a function which will play a fundamental role in our analysis, and for which geometric convergence of its factorial series comes ``natively''. We use the following representation of $\Phi$, see \cite{nist} 25.14.5:

\begin{equation}\label{ClassicalLer}
    \Phi(z,s,x)=\frac{1}{\Gamma(s)}\int_0^\infty \frac{ p^{s-1}e^{-xp}}{1-ze^{-p}}d p, \qquad    \Re s>0, \quad \Re x>0,\quad  z\in \CC\setminus[1,\infty)
\end{equation}
 For our purposes we are interested in fixing the second parameter $s=1$ and once again we use the change of variables $u={\rm e}^{- p}$ to obtain 
\begin{equation}
   \label{eq:LerInt} \Phi(z,1,x)= \int_0^1 \frac{ u^{x-1}}{1-zu}du
\end{equation}

 From \eqref{eq:fct2} we have for $z\in\DD$ and $x\notin \ZZ\setminus\NN$ we have, as $n\to \infty$, 
   \begin{equation}
   \label{PhiHorn}
     \Phi\left( \frac{z}{z-1},1,x\right)=(1-z)\sum_{k=0}^{n}z^k\frac{k!}{(x)_{k+1}}+z^nn^{-x}\Gamma(x)(z+o(1))
   \end{equation}
(for the proof see Lemma\,\ref{RemLer}). 
We see that for $|z|<1$, the domain of convergence in the $x$-plane contains the closed disk of radius one centered at 1. Thus, the Lerch function has a geometrically convergent classical factorial expansion.  It is from this object that we build our expansions which can handle functions that factorial expansions could not.

\section{Dyadic decompositions: achieving geometric convergence and extending the domain of convergence} \label{S3}

We now ask the following question; does there exist a means of improving the domain of convergence for a given classical factorial expansion?  With the use of the remarkable identity \eqref{eq:deca1},
which appears to be new, we answer this question in the affirmative.  From here, we can readily develop highly efficient methods for approximating classical special functions, and from there, for much more general ones.  

\subsection{Dyadic decomposition of the Cauchy kernel}

\begin{Lemma}[{\em Dyadic identity}]\label{L1}{\rm 
 The following identity holds in $\CC$:
\begin{equation}
      \label{eq:deca1}
      \frac1p-\left({\frac {1}{1-\mathrm{e}^{-p}}-\sum_{k=1}^\infty{\frac{2^{-k}} {1+\mathrm{e}^{-p/2^k} }}}\right)=0
    \end{equation}
   as the left hand side in \eqref{eq:deca1} has only removable singularities.
   
  Furthermore,
  
  \begin{itemize}
      \item  if we truncate \eqref{eq:deca1} then we have  
      \begin{equation}
  \label{eq:deca1traw}
   \frac1p=\frac {1}{1-\mathrm{e}^{-p}}-\sum_{k=1}^n{\frac{2^{-k}} {1+\mathrm{e}^{-p/2^k} }}+   \rho_{n+1}
   (p)  
   \end{equation}
 where 
  
      \begin{equation}
  \label{eq:deca1t}
    \rho_{n+1}(p)=\frac 1{2^n}\left( \frac 1{p/2^n} -\frac 1{1-\mathrm{e}^{-{p}/{2^n}}}\right)
    \end{equation} 
    as an equality of meromorphic functions.
   
   \item For any compact set $K\subset\CC$, if $n$ is large enough such that if $p\in K$ then $|p|<2^{n+1}\pi$,
 $\rho_{n+1}(p)$ is analytic in $K$ and uniformly bounded:  
 \begin{equation}\label{estrhoK}
 | \rho_{n+1}
 (p)|=2^{-n-1}\left(1+o(1)\right)\ \ (n\to \infty),\ \ \ p\in K
 \end{equation}

\item    The remainder $\rho_{n+1}
(p)$ also satisfies: 
 \begin{equation}
 \label{unif:estim:rhon}
 | \rho_{n+1}
 (p)|\le a_0\, 2^{-n}\ \ \ \ \text{for }|p|\le 2^n
 \end{equation}
 where
  \begin{equation}\label{valc1}
 a_0=\max_{|q|\le 1}\ \left| \frac 1q -\frac{1}{1-\mathrm{e}^{-q}}\right|
 \end{equation}
 \end{itemize}
 }
 \end{Lemma}

{\em Proof of Lemma\,\ref{L1}.}

 The proof is elementary:
 \begin{equation}
   \label{eq:decx2} \frac{1}{1-x}=\frac{2}{1-x^2}-\frac{1}{x+1}=\frac{4}{1-x^{4}}-\frac{2}{x^2+1}-\frac{1}{x+1}=\ldots=\frac{2^n}{1-x^{2^n}} -\sum_{j=0}^{n-1}\frac{2^j}{1+x^{2^j}} 
 \end{equation}
which implies, with $x=\mathrm{e}^{-p/2^n}$,
\begin{equation}
   \label{eq:decx3} 
   \frac{2^{-n}}{1-\mathrm{e}^{-\frac{p}{2^n}}}=\frac{1}{1-\mathrm{e}^{-p}}-\sum_{k=1}^n{\frac{2^{-k}} {\mathrm{e}^{-\frac{p}{2^k}}+1 }}
 \end{equation}
which implies \eqref{eq:deca1t}. From \eqref{eq:deca1t} and \eqref{eq:decx3} we see that $\rho_n(p)$ is analytic for $|p|<2^{n+1}\pi$ ($p=0$ is a removable singularity of $\rho_n$) and for $n$ large enough (so that this disk contains $K$). The bound \eqref{estrhoK} for $| \rho_n(p)|$ is immediate from \eqref{eq:deca1t}.

 To estimate \eqref{unif:estim:rhon} we take \eqref{valc1} in \eqref{eq:deca1t}.
$\Box$

In Corollary\,\ref{dyadic Cauchy} we shift $p$, to obtain the dyadic decomposition of the Cauchy kernel. It is also useful in applications to rotate $p$, so that the poles of the denominator, which in  \eqref{eq:deca1traw} are along $i\RR$, can be placed along another line. 
Let then $\beta\ne 0$. The linear affine transformation  $p\to \beta p-\beta s$ gives the following generalization of Lemma\,\ref{L1} for the Cauchy kernel.

\begin{Corollary}[{\em Dyadic decomposition of the Cauchy kernel}]\label{dyadic Cauchy}
{\rm Assume $\beta\ne 0$. 

The following identity holds in $\CC^2$:
\begin{equation}
  \label{eq:deca2}
 \frac{1}{s-p}-\left( {-\frac{\beta  \mathrm{e}^{-\beta s}}{\mathrm{e}^{- \beta s }-\mathrm{e}^{-\beta p}   }  }+ {\sum _{k=1}^{\infty } \frac{2^{-k} \beta \mathrm{e}^{- \beta s/2^{k} }}{\mathrm{e}^{-
   \beta s/2^{k}}+\mathrm{e}^{- \beta  p/2^{k}}}}\right)=0
\end{equation}
 as the left hand side in \eqref{eq:deca2} has only removable singularities.
 
 We have
   \begin{equation}
  \label{eq:deca1trawbeta}
 \frac{1}{s-p}= {-\frac{\beta  \mathrm{e}^{-\beta s}}{\mathrm{e}^{- \beta s }-\mathrm{e}^{-\beta p}   }  }+ {\sum _{k=1}^{n} \frac{2^{-k} \beta \mathrm{e}^{- \beta s/2^{k} }}{\mathrm{e}^{-
   \beta s/2^{k}}+\mathrm{e}^{- \beta  p/2^{k}}}}+ \rho_{n+1}(p,s;\beta)
   \end{equation}
 where
 \begin{equation}
  \label{eq:remnbeta}
    \rho_{n+1}(p,s;\beta
    )=\frac{\beta}{2^n}\left( \frac 1{\beta(s-p)/2^n} +\frac{  \mathrm{e}^{-{\beta s}/{2^n}}}{\mathrm{e}^{-\beta {s}/{2^n}}-\mathrm{e}^{-\beta{p}/{2^n}}}  \right)
    \end{equation}
       as an equality of meromorphic functions.

Also, $\rho_{n+1}(p,s;\beta)$ is analytic in both $p$ and $s$ throughout any compact set $K$ 
if $n$ is large enough, and it is uniformly bounded:  
\begin{equation}\label{estrhoKbeta}
| \rho_{n+1}(p,s;\beta
)|=|\beta|2^{-n-1}(1+o(1))\ \ (n\to \infty)\ \ \ \ s, p\in K
\end{equation}
  Moreover, we have
    \begin{equation}
  \label{estimCor5}
  | \rho_{n+1}(p,s;\beta
  )|\leq |\beta|\,a_0\, 2^{-n}\ \ \ \text{for }  |p-s|<2^n/|\beta|
    \end{equation}
    }
\end{Corollary}

\begin{proof} This is an immediate calculation, by replacing $p$ with $\beta p-\beta s$ in Lemma\,\ref{L1}. 

Noting that $\rho_{n}(p,s;\beta
)=-\beta\rho_n\left(\beta (p-s)\right)$ the estimates for the remainder follow.
\end{proof}

\begin{Remark}
\label{Rem2}   
{\rm    For $p\ne s$, the denominators in \eqref{eq:deca1trawbeta} vanish only for $p$ along a line. Varying the parameters (the slope of $\beta$) we can arrange that the denominators do not vanish in the complex plane cut along a ray of our choosing, and moreover, their absolute values are bounded below by a positive constant.  Then, if $p$ is such that for some constant $c_0>0$ we have
   \begin{equation}
   \label{convassump}
 |1+ \mathrm{e}^{ \beta s/2^{k}}\mathrm{e}^{- \beta  p/2^{k}}|>c_0>0,\ k=1,2,\ \ \ \ldots\text{ and }\ \  \ |1-\mathrm{e}^{ \beta s }\mathrm{e}^{-\beta p}|>c_0
 \end{equation}
 then
  \begin{equation}
  \label{goodrhon}
   | \rho_{n+1}(p,s;\beta
   )|=\left|  {\sum _{k=n+1}^{\infty } \frac{2^{-k} \beta \mathrm{e}^{- \beta s/2^{k} }}{\mathrm{e}^{- \beta s/2^{k}}+\mathrm{e}^{- \beta  p/2^{k}}}}\right|\le 2^{-n}\frac{|\beta|}{c_0}
       \end{equation}
and the dyadic series 
 \begin{equation}
   \label{convdyadic}
   \frac{1}{s-p}= {-\frac{\beta  \mathrm{e}^{-\beta s}}{\mathrm{e}^{- \beta s }-\mathrm{e}^{-\beta p}   }  }+ {\sum _{k=1}^{\infty } \frac{2^{-k} \beta \mathrm{e}^{- \beta s/2^{k} }}{\mathrm{e}^{- \beta s/2^{k}}+\mathrm{e}^{- \beta  p/2^{k}}}}
    \end{equation}
 converges geometrically.
}
   \end{Remark}

\subsubsection{More about the Lerch $\Phi$ trancendent and formula \eqref{PhiHorn}}
The proof of Theorem\,\ref{ElementDecomp} is illuminated, and simplified, by formula \eqref{PhiHorn} which we state in detail here and prove in \S\ref{ProofL11}.

Denote $\ZZ_-=\{0,-1,-2,\ldots,\}$.

\begin{Lemma}\label{RemLer}
   For $|z|<1$ and $x\in\CC\setminus\ZZ_-$ we have 
   \begin{equation}
   \label{factorialseriesforPhi}
    \Phi\left( \frac{z}{z-1},1,x\right)=(1-z)\sum_{j\ge 0}z^j\frac{j!}{(x)_{j+1}}
 \end{equation}
and the series converges absolutely.

The remainder 
 \begin{equation}
   \label{rhon}
   \rho_{n+1,0}
   (z,x):=(1-z)\sum_{k= n+1}^\infty z^k\frac{k!}{(x)_{k+1}}
   \end{equation}
   satisfies the following estimates:
   
   (i) for large enough $n$,
      \begin{equation}
\label{rhonestimate}
\rho_{n+1,0}
(z,x) =z^n\frac{n!}{(x)_{n+1}}\left(z-u_n\right),\ \ \  \text{with }|u_n|< \frac 1n \frac1{1-\lambda} |xz|\, \sup _{n> M} \frac{n-1}{|x+n|},\ \ \ \ \ \ \text{for } n\ge M
\end{equation}
where $M$ is determined as follows. Let $x\not\in \ZZ_-$, $x=|x|e^{i\alpha}$, and $\lambda$ is such that $|z|\le\lambda<1$.

\ \ \ \ \ \ \ (i1) $M$ is such that 
\begin{equation}
 \label{valM}
\text{ if }\ \ \cos^2\alpha-1+|z|^2/\lambda^2\ge 0\ \text{ and }\ \cos\alpha<0
\ \ \ \text{then }\ \ M>\frac {|x|}{-\cos\alpha -\sqrt{\cos^2\alpha-1+|z|^2/\lambda^2}}\ \ \ \ \end{equation}

\ \ \ \  \ \ \ (i2) If the assumption in \eqref{valM} does not hold, then we let $M=1$.
 
 (ii) for any $n$, the remainder has the expression
   \begin{equation}
\label{newrhonestimate}
\rho_{n+1,0}
(z,x) =z^{n+1}\frac{(n+1)!}{(x)_{n+1}}(1-z)\, \int_\mathcal{C}\frac{ {\rm{e}}^{-px}}{{\rm{e}}^{p(n+1)}(1-z+z{\rm{e}}^{-p})^{n+2}}\, dp \ \ \ \ \ \text{for }x\in\CC\setminus(-\infty,0]
\end{equation}
  where $\mathcal{C}=[0,+\infty)$ if $\Re (x)>0$, while if $\Re (x)\le 0$ then $\mathcal{C}$ is the segment $[0,p_z] $ followed by $p_z+{\rm{e}}^{i\phi}\RR_+$ for a $p_z$ satisfying $p_z\ge 0$, $p_z>-\ln\left({|z|}^{-1}-1\right)$ and $\phi\in(-\tfrac{\pi}2,\tfrac{\pi}2)$ is such that $\Re (x{\rm{e}}^{i\phi})>0$.
 
(iii) for large $n$ and $|z|<1$ we have 
 \begin{equation}
\label{newrhonestimate2}
\rho_{n+1,0}
(z,x) \sim z^{n+1}(n+1)^{-x}\Gamma(x)
\end{equation}

\end{Lemma}
 
 \begin{Note} {\rm If $x$ is in the left half-plane, we note that \eqref{newrhonestimate2} 
 requires that $n$ be large relative to $x$. For smaller $n$, the remainder is estimated using formula \eqref{newrhonestimate}.}
 
 \end{Note}
 
 The proof of Lemma\,\ref{RemLer} is given in \S\ref{ProofL11}.

\subsection{Dyadic expansions in a cut plane for {\em function elements}}
 \subsubsection{Function element}\label{FuncEl}
In \S\ref{Resfun} we show that very general classes of functions can be decomposed in terms of simpler functions, namely functions $F$ with only one singularity on the first Riemann sheet in the Borel plane, at say $p=p_0$, such that $(p-p_0)^mF(p)$ is locally bounded (for some $m$). 

\bd
    A function $F$ is called a {\em function element} if \newline
   (i) $F$ is analytic at the origin and a cut plane, i.e. analytic in a domain of the form $\mathscr{D}=\CC\setminus l_\omega$ where $l_\omega$ is a half-line originating at $\omega\in\CC\setminus\{0\}$, and \newline
   (ii) $F$ decays in $\mathscr{D}$ as $|p|\to \infty$.\footnote{For functions occuring in many applications such as generic ODEs and difference equations, the function elements are in fact {\em resurgent} but we do not analyze these further features here.}
\ed

More precise statements are made when needed for proofs, see for example the hypothesis of Theorem\,\ref{ElementDecomp}.

\subsubsection{Dyadic series for function elements}

Theorem\,\ref{ElementDecomp} finds the dyadic series for functions which are Laplace transforms of function elements. It also shows that the dyadic series converges in a cut plane. The theorem also estimates the errors when a truncation is used to approximate the function.

Theorem\,\ref{ElementDecomp}  treats the case when the singularity of $F$ is a branch point. The special case when the singular point is a pole is simpler, and we illustrate its treatment in sections \S\ref{EiStokes} and \S\ref{Eiaway}.
In any case, integration by parts transforms such a function to a function with (logarithmic) branch point.

To motivate the setting of Theorem\,\ref{ElementDecomp} and explain how it can be used for various function elements we note the following. First, by changes of variable the singularity can be placed anywhere in the complex plane.  
Now suppose $F(p)$ has an (integrable) branch point singularity at $p=1$ and $F$ decays at $\infty$ fast enough to be $L^1$ along rays. 
Consider its Laplace transform, $f(x)$ defined for $x\in\CC$ with $\arg x=0+$ by the Laplace transform of $F$ along $e^{i\phi}\RR_+$ where $\phi=0-$:
$f(x)=\int_0^{\infty e^{i\phi}}\, e^{-xp}F(p)\, dp$.  For other values of $x$, $f(x)$ is defined by analytic continuation.
Since $F$ is assumed to decay at infinity, the analytic continuation of $f(x)$ for larger arguments of $x$ can be obtained from the Laplace integral by simultaneously rotating $\phi$ clockwise (as long as ${\rm{e}}^{-i\phi}$ does not cross $\RR_+$)
 and $x$ anticlockwise in such a way that $\Re(px)>0$ throughout the rotation.
This is accomplished 
 by ensuring $\arg x\in( -\phi-\tfrac\pi 2, -\phi+\tfrac\pi 2)$.

This motivates the definition of $f$, the Laplace transform of $F$, as given by
\begin{equation}
 \label{Laplaceb}
  f(x)=\int_0^{\infty {\rm{e}}^{-ib}}\, {\rm{e}}^{-xp}F(p)\, dp\ \ \ \ \ \ \text{ for }\arg x \in\left( b-\frac\pi 2, b+\frac\pi 2\right)
  \end{equation}

\begin{Theorem}

\label{ElementDecomp}

Let $\beta\ne 0$ so that $\beta=|\beta| {\rm{e}}^{ib}$ with $b\in\left[ \tfrac{\pi}2,\tfrac{3\pi}2\right]$  and
\begin{equation}\label{condbeta}
{\rm{e}}^{-\Re\beta}-2\cos(\Im\beta):=c_\beta>0,\ \ \ \ \ |\beta|\le\pi
\end{equation}
(the limitation on $|\beta|$ is for convenience). 

 Let $\theta$ be the angle in the right half-plane so that $b+\theta=\pi(\text{mod } 2\pi)$.
 
 Assume that a  function $F$ in the Borel plane has the following properties: 
 \begin{enumerate}
      \item $F$ has exactly one singularity: an integrable branch point (we place it conventionally at $p=1$).
      \item $F$ decays at $\infty$: $|F(p)|\le C|p|^{-\alpha}$ (with $\alpha>1$) for large $|p|$, and $F$ is $L^1_{\rm loc}$.
      \item $F$ is analytic in the cut plane $\CC\setminus [1+{\rm{e}}^{i\theta}[0,\infty) ]$, and can be analytically continued through both sides of the cut.
           \end{enumerate}

Then the function $f(x)$ defined by \eqref{Laplaceb}
 has the dyadic expansion, for all $x\in\CC\setminus {\rm{e}}^{ib}(-\infty,0]$,
\begin{equation}
\label{dyadicf}
f(x)=  \sum_{m=1}^{n-1}\frac{(m-1)!}{(x/\beta)_{m}}\, d_{m,0}+\sum_{k=1}^{N-1} \sum_{m=1}^{\ell}\frac{(m-1)!}{(2^kx/\beta)_{m}}\, d_{m,k}+\mathcal{R}_{n,N,\ell}\left(x,\beta\right)
\end{equation}
where, denoting by $\Delta F(1+t {\rm{e}}^{i\theta})$ the branch jump of $F$, $F(1+t {\rm{e}}^{i\theta^+})-F(1+t {\rm{e}}^{i\theta^-})$ and $s=1+t {\rm{e}}^{i\theta}$, the coefficients of the series have the expressions
\begin{equation}
\label{dmk}
d_{m,0}=\frac {{\rm{e}}^{i\theta}}{2\pi i}\int_0^\infty \, \Delta F(1+t {\rm{e}}^{i\theta})\, \frac{{\rm{e}}^{\beta s (m-1)}}{({\rm{e}}^{\beta s}-1)^m} dt,\ \ \ \ \ \ d_{m,k}=\frac {{\rm{e}}^{i\theta}}{2\pi i}\int_0^\infty\,  \Delta F(1+t {\rm{e}}^{i\theta})\,  \frac{{\rm{e}}^{\beta s (m-1)/2^k}}{({\rm{e}}^{\beta s/2^k}+1)^m}dt 
\end{equation}
and $\mathcal{R}_{n,N,\ell}\left(x,\beta\right)$
has the expression 
\begin{equation}
\label{RnNell}
\mathcal{R}_{n,N,\ell}\left(x, \beta\right)
=\frac {{\rm{e}}^{i\theta}}{2\pi i}\int_0^\infty\,  \Delta F(1+t {\rm{e}}^{i\theta})\,\left( -\rho_{n,0}(t, x;\beta
) +\sum_{k=1}^{N-1} \rho_{\ell,k}(t, x; \beta
)+R_N(t, x; \beta
)\right)dt
\end{equation}
where the remainder terms seen in \eqref{RnNell} are defined by

\begin{enumerate}
\item 
 \begin{equation}
\label{newrhonestimateThm}
\rho_{n,0}
(t,x; \beta) =\frac{(-1)^n {\rm{e}}^{\beta sn}n!}{(x/\beta)_{n}}\, \int_{\Gamma_{c}}\frac{ {\rm{e}}^{-q(x/\beta+n)}}{(1-{\rm{e}}^{\beta s}{\rm{e}}^{-q})^{n+1}}\, dq \ \ \ \ \ \text{for } n=1,2,\ldots
\end{equation}
\item    
\begin{multline}
\label{newrhonestimateThmL}
    \rho_{\ell,k} (t,x; \beta
)= {\rm{e}}^{\beta s\ell /2^k}  \frac{\ell!}{(2^kx/\beta)_{\ell}}\, \int_{\Gamma_{c}} \frac{ {\rm{e}}^{-q(2^kx/\beta+\ell)}}{\left({\rm{e}}^{\beta s/2^k}{\rm{e}}^{-q}+1\right)^{\ell+1}}\, dq \ \ \ \text{for} \,\ell\in \NN, \ k=1,\ldots N-1
\end{multline}
see \eqref{eq:remnbeta}

\item 
 
For $x\notin (-\infty,0]\beta$ we have
\begin{equation}
\label{RNEstThm}
    R_N(t,x; \beta
) =\frac{1}{\beta}\int_{\Gamma_{c}} \mathrm{e}^{-xq/\beta} \rho_{N}(q/\beta,1+t {\rm{e}}^{i\theta};\beta)dq\\\ 
\end{equation}
The contour  \( \Gamma_{c} \) is chosen so that along it we have: (a)  $\Re(xq/\beta)>0$ for large $q$ and (b) the function $\rho_N$ (defined by \eqref{eq:remnbeta} with $1+t {\rm{e}}^{i\theta}=s$) is analytic. More precisely,
given \( x\in \mathbb{C} \setminus \beta(-\infty,0] \) and any  
\( c\in \left(0,\min\left\{1,\frac{1}{2}\operatorname{dist}(x,\beta(-\infty,0])\right\}\right) \),  define

\begin{equation}
\label{RemContGam}
    \Gamma_{c}=\begin{cases}
        \RR^+ & \Re(x/\beta)>c\\[6pt]
       \mathrm{e}^{-i\pi/2}\RR^+ &   \Re(x\mathrm{e}^{-i\pi/2}/\beta)>c\\[6pt]
        \mathrm{e}^{i\pi/2}\RR^+ &   \Re(x\mathrm{e}^{i\pi/2}/\beta)>c
        
    \end{cases}
\end{equation}  
\end{enumerate}
(see also Remark\, \ref{Remark10}).

Below, $s=1+t{\rm{e}}^{i\theta}$.  The remainders satisfy, for all $t\ge 0$: 

\begin{itemize}

 \item for large $\ell$ and $k=1,2,\ldots$ 
 
\begin{equation}
        \label{rholargel}
\rho_{\ell,k}(t, x;\beta
)=\frac{{\rm{e}}^{\beta s \ell/2^k}}{({\rm{e}}^{\beta s/2^k}+1)^{\ell}} \frac{(\ell-1)!}{(2^kx/\beta)_\ell}\left(1+o(1)\right)\;\text{for }x\in \Omega_c 
\end{equation}

\item while for $k=0$ and large $n$, 
\begin{equation}
\label{rholargen}
\rho_{n,0}(t, x;\beta
)=\frac{{\rm{e}}^{\beta s n}}{({\rm{e}}^{\beta s}-1)^{n}} \frac{(n-1)!}{(x/\beta)_n}\left(1+o(1)\right)\; \text{for }\;x\in\Omega_c
\end{equation}

\item Also, for large $N$, $x\in\Omega_c$ and $c>0$,
\begin{equation}
\label{eq:singrem}
    |R_N(t, x;\beta
    )|\le
     \begin{cases}
      \frac1{c_02^{N-1}\Re (x/\beta)}& \Re(x/\beta)>c\\[10pt]
    \frac{1}{c_02^{N-1}\Re(x\mathrm{e}^{-i\pi/2}/\beta)}& x\in\{\Re(x\mathrm{e}^{-i\pi/2}/\beta)> c\}\cap\{\Re (x/\beta)\leq c\}\\[10pt]
    \frac{1}{c_02^{N-1}\Re(x\mathrm{e}^{i\pi/2}/\beta)}& x\in\{\Re(x\mathrm{e}^{i\pi/2}/\beta)> c\}\cap\{\Re (x/\beta)\leq c\}
    \end{cases}
 \end{equation} 
\item For any $n$ and $\ell$ (not necessarily large), $x\in\Omega_c$ and $c>0$ the remainders satisfy
\begin{equation}
\label{finalestimrhoeln} 
|\rho_{n,0}(t,x;\beta)|
\le 
\begin{cases}
\frac{c_1^{n}n!}{c_0 |(x/\beta)_{n}|\Re({ x/\beta})} & \Re \left(x/\beta\right)>c
\\[10pt]
\frac{c_1^{n}n!}{c_0|(x/\beta)_{n}|\Re(x\mathrm{e}^{-i\pi/2}/\beta)} & x\in\{\Re(x\mathrm{e}^{-i\pi/2}/\beta)> c\}\cap\{\Re (x/\beta)\leq c\}\\[10pt]
\frac{c_1^{n}n!}{c_0|(x/\beta)_{n}|\Re(x\mathrm{e}^{i\pi/2}/\beta)}
& x\in\{\Re(x\mathrm{e}^{i\pi/2}/\beta)> c\}\cap\{\Re (x/\beta)\leq c\}
\end{cases}
\end{equation}

respectively
\begin{equation}
\label{finalestimrhoelk} 
|\rho_{\ell,k}(t,x;\beta)|\le 
\begin{cases}
\frac{c_1^{\ell}\ell !}{c_02^k{|(2^kx/\beta)_{\ell}|\Re({ x/\beta})}}   & \Re\left(x/\beta\right)>c, \;k\in\ZZ^+
\\[10pt]
\frac{c_1^{\ell}\ell !}{c_0 2^k{|(2^kx/\beta)_{\ell}|}\Re(x\mathrm{e}^{-i\pi/2}/\beta)}\ & x\in\{\Re(x\mathrm{e}^{-i\pi/2}/\beta)> c\}\cap\{\Re (x/\beta)\leq c\}, \;k\in\ZZ^+\\[10pt]
\frac{c_1^{\ell}\ell !}{c_0 2^k{|(2^kx/\beta)_{\ell}|}\Re(x\mathrm{e}^{i\pi/2}/\beta)}& x\in\{\Re(x\mathrm{e}^{i\pi/2}/\beta)> c\}\cap\{\Re (x/\beta)\leq c\}, \;k\in\ZZ^+
\end{cases}
\end{equation}
(see also Remark\,\ref{RemIneqRem}). 
\end{itemize}

Letting $n,N,\ell\to \infty$ we see that the series \eqref{dyadicf} converges absolutely.
\end{Theorem}

The proof is found in \S\ref{PfTh8}.

\begin{Remark}\label{Remark10}
The choice of $\Gamma_c$ ensures that the integral representations \eqref{newrhonestimateThm}, \eqref{newrhonestimateThmL}, and \eqref{RNEstThm} are holomorphic in a half-plane neighborhood of \( x \), which  may then be analytically continued by contour deformation to the region \( \Omega_c \)  consisting of  the union of three half planes:  
     \begin{equation}\label{DomHolRem}
        \Omega_c=\{\Re(x/\beta)>c\}\cup \{\Re(x\mathrm{e}^{-i\pi/2}/\beta)>c\}\cup \{\Re(x\mathrm{e}^{i\pi/2}/\beta)>c\}
    \end{equation}
    This region is the complement of a closed rectangular strip centered on the ray $\beta(-\infty,0]$ whose boundary is at a positive  distance from the ray.  It is this positive separation from the cut that ensures our remainder estimates hold on compact subsets of the domain.  
\end{Remark}

\begin{Remark}\label{RemIneqRem}
    Since $c>0$ is arbitrary strict inequalities can be replaced by non-strict ones.  
\end{Remark}

\begin{Note}
 In practice it is useful to take the order $\ell$ of truncation in \eqref{dyadicf} to be $k$-dependent: when the goal is to achieve a desired precision in approximating $f(x)$, the larger $k$ is, the smaller $\ell$ can be. See an example in Fig.\,\ref{fig12}.
 \end{Note}
 See also Note\,\ref{Note12} for the relation between $x$ and the number of terms needed.

\begin{Theorem}\label{choosebeta}
 The poles of the dyadic series \eqref{dyadicf} are for $x$ on the ray $(-\infty,0]\beta$.  Convergence is uniform on compact sets contained in the complement of the ray. 
 
\end{Theorem}

\begin{Note}\label{choosebeta11} 
 By choosing $\beta$ we can choose the cut plane where we obtain the approximation of $f(x)$. Hence the domain of convergence of the dyadic series is the complex plane without a cut that can be placed anywhere in the closed right half-plane by an appropriate choice of $\beta$.
 \end{Note}

The proof of Theorem\,\ref{choosebeta} is found in \S\ref{PfofThchoosebeta}.

\subsection{Proof of Theorem\,\ref{ElementDecomp}} \label{PfTh8}

  Note that the path of integration in \eqref{Laplaceb}, $p=|p|{\rm{e}}^{-ib}$ and the cut $s=1+t{\rm{e}}^{i\theta}$, $t\ge 0$ do not intersect. Indeed, noting that our choice of $\theta$ implies that $\sin\theta=\sin b$ and $\cos\theta=-\cos b$, a point on the intersection would satisfy $|p|+t={\rm{e}}^{ib}$ which is not possible for our restriction on $b$.

Let $\mathcal{C}_1$ be a simple closed contour in the cut plane $\CC\setminus[1+ {\rm{e}}^{i\theta}[0,\infty) ]$ and $p$ a point inside $\mathcal{C}_1$.
Using the Cauchy formula, then deforming the path of integration to a Hankel contour hanging around the cut we have
\begin{equation}
\label{formulaFCauchy}
F(p)=\frac 1{2\pi i}\oint_{\mathcal{C}_1}\frac {F(s)}{s-p}ds=\frac 1{2\pi i}\int_1^{1+\infty {\rm{e}}^{i\theta}}\frac {F(s+)-F(s-)}{s-p}ds=\frac {{\rm{e}}^{i\theta}}{2\pi i}\int_0^\infty\,  \frac {\Delta F(1+t {\rm{e}}^{i\theta})}{1+t {\rm{e}}^{i\theta}-p} dt
\end{equation}
As noted, $1+t {\rm{e}}^{i\theta}-p\ne 0$ for $p\in  {\rm{e}}^{-ib}[0,+\infty)$ and $t\geq 0$ .

Taking the Laplace transform \eqref{Laplaceb} in \eqref{formulaFCauchy}, then interchanging the order of integration we obtain 

\begin{equation}
\label{formulaf}
f(x)= \frac {{\rm{e}}^{i\theta}}{2\pi i}\int_0^\infty\,   \Delta F(1+t {\rm{e}}^{i\theta})dt\int_0^{\infty {\rm{e}}^{-ib}}\ \frac{{\rm{e}}^{-xp}} {1+t {\rm{e}}^{i\theta}-p} dp \ \ \ \text{for}\quad \arg x\in \left(b-\frac{\pi}{2},b+\frac{\pi}{2}\right)
\end{equation}

In the integral above denote  
\begin{equation}\label{defJdet}
J(t,x):=\int_0^{\infty {\rm{e}}^{-ib}}\frac{{\rm{e}}^{-xp}} {1+t {\rm{e}}^{i\theta}-p}dp
\end{equation}
where we now use the dyadic series \eqref{eq:deca1trawbeta} for $s=1+t {\rm{e}}^{i\theta}$. This series converges geometrically. Indeed, along the path of integration $p\in {\rm{e}}^{-ib}[0,+\infty)$, and for $s=1+t{\rm{e}}^{i\theta}\ (t\geq 0)$ we have $1-{\rm{e}}^{\beta s}{\rm{e}}^{-\beta p}\neq0$ and $|1+e^{\beta s/2^k}e^{-\beta p/2^k}|>1$ for all $k\geq 1$ (see also Remark\,\ref{Rem2}).

 Therefore we can interchange the order of summation with integration. Changing the variable of integration $p=q/\beta$ and keeping track of the remainder, we have
\begin{equation}
\label{dyaJ}
J(t,x)=   {- \int_0^{\infty}  \mathrm{e}^{-xq/\beta}\frac{ 1}{1-\mathrm{e}^{ \beta (1+t e^{i\theta}) }\mathrm{e}^{-q}   }  }dq+ {\sum _{k=1}^{N-1 }  \int_0^{\infty}  \mathrm{e}^{-xq/\beta}\frac{2^{-k} }{1+\mathrm{e}^{\beta (1+t  \mathrm{e}^{i\theta})/2^{k}}\mathrm{e}^{- q/2^{k}}}} dq+R_N(t,x;\beta
)
\end{equation}
where
\begin{equation}
\label{Rn}
R_N(t,x;\beta
) =\frac{1}{\beta}\int_0^{\infty}  \mathrm{e}^{-xq/\beta} \rho_{N}(q/\beta,1+t e^{i\theta};\beta
)dq\ \ \ \ \ \ \text{ for }\arg x \in\left( b-\frac\pi 2, b+\frac\pi 2\right)
\end{equation}
with $\rho_{N}(p,s;\beta)
$ as in Corollary\,\ref{dyadic Cauchy} (for $n+1=N$).  To extend the domain of \eqref{Rn} we use contour deformation as in \eqref{RNEstThm}.

After further changing the integration variable $q/2^k$ to $q$ in \eqref{dyaJ} we obtain
\begin{equation}
\label{dyaJac}
J(t,x)=   {- \int_0^{\infty}  \mathrm{e}^{-xq/\beta}\frac{ 1}{1-\mathrm{e}^{ \beta (1+t e^{i\theta}) }\mathrm{e}^{-q}   }  }dq+ {\sum _{k=1}^{N-1 }  \int_0^{\infty}  \mathrm{e}^{-2^kxq/\beta}\frac{1}{1+\mathrm{e}^{\beta (1+t  \mathrm{e}^{i\theta})/2^{k}}\mathrm{e}^{- q}}}dq+R_N(t,x;\beta
)
\end{equation}

In \eqref{dyaJac} we now use the integral representation \eqref{IntPhi} of the Lerch $\Phi$ transcendent and obtain
\begin{equation}
\label{serJPhi}
J(t,x)= -\Phi\left(e^{\beta (1+t e^{i\theta})}, 1, \frac x\beta\right)+\sum_{k=1}^{N-1} \Phi\left( -e^{\beta (1+t e^{i\theta})/2^k}, 1, \frac {2^kx}\beta\right)+R_N(t,x;\beta
)
\end{equation}

Note that due to \eqref{condbeta} we have
\begin{equation}
\label{estimz0}
\left| \tfrac{1}{1-e^qe^{-\beta (1+te^{i\theta})}} \right| \le \frac1{\sqrt{1+c_\beta}}<1,\ \ \text{ for all }q,t\ge0
\end{equation}
 Indeed, $| 1-e^qe^{-\beta (1+te^{i\theta})}|^2=1+\rho(\rho-2\cos\alpha)$ where $\alpha=|\beta|\sin b$ and $\rho=e^{q+|\beta|(t-\cos b)}\ge e^{-|\beta|\cos b}= 2\cos\alpha+c_\beta$ (and, of course, $\rho\ge 1$).

Using \eqref{factorialseriesforPhi} for $z=\tfrac{1}{1-e^{-\beta (1+te^{i\theta})}}$ (we have $|z|<1$ from \eqref{estimz0}) and, with the notation $s=1+t e^{i\theta},\ t\ge 0$ we obtain, after truncating the series,
\begin{multline}
\label{serP0}
\Phi\left( e^{\beta s},1,x/\beta\right)=\frac{-1}{e^{\beta s}-1}\sum_{j=0}^{n-1}\frac{e^{\beta s j}}{(e^{\beta s}-1)^j}\frac{j!}{(x/\beta)_{j+1}}+ \rho_{n,0}(t,x;\beta
),\ \ \ (s=1+t e^{i\theta},\ t\ge 0)\\ 
\text{where}\ \ \rho_{n,0}(t,x;\beta
)=\frac{e^{\beta s (n-1)}}{(e^{\beta s}-1)^{n-1}} \frac{(n-1)!}{(x/\beta)_n}\,\left(\frac{e^{\beta s}}{e^{\beta s}-1} +e_{n,0}\right)
\end{multline}
where, for large $n$, $|e_{n,0}|$ satisfies estimates similar to those of $u_n$ in \eqref{rhonestimate}, for all $t\ge 0$.

For moderate $n$, the integral representation for the remainder \eqref{newrhonestimateThm} is obtained by applying \eqref{newrhonestimate} to \eqref{serP0}.  Unlike in \eqref{newrhonestimate} we use an alternative contour deformation defined by the family of contours $\Gamma_c$ (see \eqref{RemContGam}).  We obtain,
\begin{equation}
\label{rhoelnestimate}
\rho_{n,0}(t,x; \beta
) =-
e^{\beta sn}\frac{n!}{(x/\beta)_{n}}\, \int_{\Gamma_{c}} \frac{ e^{-q(x/\beta+n)}}{\left(e^{\beta s}e^{-q}-1\right)^{n+1}}\, dq,\quad n=1,2,\ldots
\end{equation}

Note that 
\begin{equation}
\label{estimzk}
\left| \frac{1}{1+e^qe^{-\beta (1+te^{i\theta})/2^k}}\right| \le\frac1{\sqrt{2}}<1\ \ \text{ for all }q\ge 0,t\ge 0,\ k\ge 1
\end{equation}

 Indeed, $\left|1+e^q e^{-\beta (1+te^{i\theta})/2^k}\right|^2
 =1+\rho_k^2+2\rho_k\cos(\alpha_k)$ where $\rho_k=e^{q+|\beta| (t-\cos b)/2^k}\ge 1$ and $|\alpha_k|=\frac{|\beta|}{2^k}|\sin b|\le \frac{|\beta|}{2}\le \tfrac\pi 2$ hence $\cos\alpha_k\ge 0$.

Using \eqref{factorialseriesforPhi} for $z=\tfrac1{1+e^{-\beta s/2^k}}$ (which, by \eqref{estimzk}, satisfies $|z|<1$) and denoting $s=1+t e^{i\theta}$, we obtain
\begin{multline}
\label{serPk}
\Phi\left( -e^{\beta s/2^k}, 1, \frac {2^kx}\beta\right)= \frac{1}{e^{\beta s/2^k}+1}\sum_{j=0}^{\ell-1}\frac{e^{\beta s j/2^k}}{(e^{\beta s/2^k}+1)^j}\frac{j!}{(2^kx/\beta)_{j+1}}+ \rho_{\ell,k}(t,x;\beta
),\\ 
\text{where}
\  \rho_{\ell,k}(t,x;\beta
)=\frac{e^{\beta s (\ell-1)/2^k}}{(e^{\beta s/2^k}+1)^{\ell-1}} \frac{(\ell-1)!}{(2^kx/\beta)_\ell}\left(\frac{e^{\beta s/2^k}}{e^{\beta s/2^k}+1}+e_{\ell,k}\right)
\end{multline}
where, for large $\ell$,  $|e_{\ell,k}|$ satisfies estimates similar to those of $u_n$ in \eqref{rhonestimate}, for all $t\ge 0$.

For moderate $\ell$, the integral representation for the remainder \eqref{newrhonestimateThmL} is obtained by applying \eqref{newrhonestimate} 
to \eqref{serPk}.  As in \eqref{rhoelnestimate}, we use an alternative contour deformation defined by the family of contours $\Gamma_c$ (see \eqref{RemContGam}) and obtain: 

      \begin{equation}
\label{rhoellestimate} 
\rho_{\ell,k} (t,x; \beta
)= e^{\beta s\ell /2^k}  \frac{\ell!}{(2^kx/\beta)_{\ell}}\, \int_{\Gamma_{c}}\frac{ e^{-q(2^kx/\beta+\ell)}}{\left(e^{\beta s/2^k}e^{-q}+1\right)^{\ell+1}}\, dq,\quad \ell=1,2,\ldots
\end{equation}
Using \eqref{serP0} and \eqref{serPk} in \eqref{serJPhi} we obtain
\begin{multline}\label{dyadicJdet}
J(t,x)=\frac{1}{e^{\beta s}-1}\sum_{j=0}^{n-1}\frac{e^{\beta s j}}{(e^{\beta s}-1)^j}\frac{j!}{(x/\beta)_{j+1}}- \rho_{n,0}(t,x;\beta
)\\
+\sum_{k=1}^{N-1}\left[ \frac{1}{e^{\beta s/2^k}+1}\sum_{j=0}^{\ell-1}\frac{e^{\beta s j/2^k}}{(e^{\beta s/2^k}+1)^j}\frac{j!}{(2^kx/\beta)_{j+1}}+ \rho_{\ell,k}(t,x;\beta
)\right]+R_N(t,x;\beta
)
\end{multline}
which introduced in \eqref{formulaf} gives
\begin{multline}
\label{dyadicfraw}
f(x)= \frac {e^{i\theta}}{2\pi i}\int_0^\infty\,  \Delta F(1+t e^{i\theta})\,\left[ \sum_{m=1}^{n}\frac{e^{\beta s (m-1)}}{(e^{\beta s}-1)^m}\frac{(m-1)!}{(x/\beta)_{m}}- \rho_{n,0}(t,x;\beta
) \right.\\
\left. +
\sum_{k=1}^{N-1}\left( \sum_{m=1}^{\ell}\frac{e^{\beta s (m-1)/2^k}}{(e^{\beta s/2^k}+1)^m}\frac{(m-1)!}{(2^kx/\beta)_{m}}+ \rho_{\ell,k}(t,x;\beta
)\right)+R_N(t,x;\beta
)\right]dt
\end{multline}
and  \eqref{dyadicf}, \eqref{dmk},\eqref{RnNell} follow. 

\

{\em Estimates of the remainders.} 

To estimate the remainders \eqref{rhoelnestimate} and \eqref{rhoellestimate} for $n$, respectively $\ell$, moderate we first note that there is a constant $c_0>0$ so that
\begin{equation}
\label{estdenom}
\left| 1-e^{\beta s}e^{-q}\right|>c_0>0,\ \ \ \ \left| 1+e^{\beta s/2^k}e^{-q}\right|>c_0>0\ \ \text{for all }q,t\ge 0,k \geq 1
\end{equation}

Indeed, $\left| 1-e^{\beta s}e^{-q}\right|^2=r^2-2r\cos\alpha+1$ where $r=e^{-q-|\beta|(t-\cos b)}\in (0,  e^{|\beta| \cos b}]\subset (0,1)$ and $\alpha=|\beta|\sin b\in[-\pi,\pi]$. 
If $|\alpha|\ge \pi/2$ then we can clearly take $c_0=1$. Otherwise, if $|\alpha|< \pi/2$, then $r^2-2r\cos\alpha+1=(r-\cos\alpha)^2+\sin^2\alpha\ge \sin^2\alpha=\sin^2(|\beta|\sin b) :=c_{0,1}^2$. We have $c_{0,1}>0$ for $\alpha=0$, while for $\alpha=0$, meaning $b=\pi$, we have $r^2-2r\cos\alpha+1=(1-r)^2\ge (1-e^{|\beta| \cos b})^2=(1-e^{-|\beta|})^2:=c_{0,0}^2>0$. We let $c_0=\min\{c_{0,1},c_{0,0}\}>0$.
The second inequality in \eqref{estdenom} holds for $c_0=1$ since $\left| 1+e^{\beta s/2^k}e^{-q}\right|^2=r_k^2+2r_k\cos\alpha_k+1$ where $r_k=e^{-q-|\beta|/2^k(t-\cos b)}\in(0,1)$ and $\alpha_k=|\beta|/2^k\sin b\in[-\pi/2,\pi/2]$ hence $\cos\alpha_k\ge 0$.

Let us unify the estimates \eqref{estimz0} and \eqref{estimzk} by writing
\begin{equation}\label{c1Const}
    c_1:=\max\left\{ \frac1{\sqrt{1+c_\beta}},\frac1{\sqrt{2}}\right\}<1
\end{equation}

Using \eqref{estdenom} and \eqref{estimz0} in \eqref{rhoelnestimate} we estimate and obtain \eqref{finalestimrhoeln} for $\Re (x/\beta)\geq c$.  Analogous estimates using the other two cases of \,$\Gamma_c$ yield the rest of \eqref{finalestimrhoeln}.  Similarly, using \eqref{estdenom} and \eqref{estimzk} in \eqref{rhoellestimate} we obtain \eqref{finalestimrhoelk} for $\Re(x/\beta)\geq c$ and perform the additional analysis as above to obtain the rest of estimate \eqref{finalestimrhoelk}.

Applying \eqref{rhonestimate} to \eqref{serP0} and \eqref{serPk} for large $n$ and large  $\ell$, we obtain the estimates \eqref{rholargen} and \eqref{rholargel}.

Finally, for large $N$, the case  $\Re(x/\beta)\geq c$ in \eqref{eq:singrem} follows from \eqref{Rn} and \eqref{goodrhon}.  For the other components of $\Omega_c$, straightforward estimates for the other two cases of $\Gamma_c$ imply \eqref{eq:singrem}. $\Box$

 \subsection{Proof of Theorem\,\ref{choosebeta}}\label{PfofThchoosebeta}
\begin{proof} 
Since $\left( 2^kx/\beta\right)_m=0$ only for $2^kx/\beta=-n$ with $n\in\{0,1,2,\ldots,m-1\}$ the terms in  the series are defined only when $x\not\in\beta(-\infty,0]$.

We perform a change of coordinates by rotating the $x$-plane. Let $y=x/\beta$ and $y=|y|e^{i\alpha}$ with $y\in\CC\setminus (-\infty,0]$, so that the denominators in \eqref{dyadicf} do not vanish.  Given a compact set $K\subset \CC \setminus(-\infty,0]$  choose $c>0$ small enough so that the set $\mathcal{S}_c$ which is the complement of the union of three closed half-planes :

$$\mathcal{S}_c=\CC\setminus\left(\{\Re y\geq c\}\cup\{\Im y\geq c\}\cup \{\Im y\leq-c\}\right)$$ 
is disjoint from $K$.  By breaking $K$, if necessary, into compact subsets with disjoint interiors,  we arrange that $K\subset \{\Re y\geq c\}$, or $K\subset \{\Im y\leq -c\}$ or $K\subset \{\Im y\geq c\}$ (see Remark\ref{RemIneqRem}).  Using the estimates \eqref{rholargel}, \eqref{rholargen},\eqref{eq:singrem} and assuming $l$, $n$ and $N$ are all sufficiently large, we obtain an upper bound on $\left|\mathcal{R}_{n,N,\ell}\left(\beta y,\beta\right)\right|$ of the form (see \eqref{RnNell},\eqref{estimz0} and \eqref{c1Const}):

\begin{equation}
\label{FullRemEst}
      \left|\mathcal{R}_{n,N,\ell}\left(\beta y,\beta\right)\right| \leq M\left(\frac{c_1^{n}(n-1)!}{{|(y)_{n}|}}   +\sum_{k=1}^{N-1} \frac{c_1^{\ell}(\ell-1)!}{{|(2^ky)_{\ell}|}}   +\frac{1}{c_0 2^{N-1}\Re(\sigma (y))}\right)
\end{equation}
throughout $\mathcal{S}_c$ for some $M>0$ independent of $n,\ell,k$ and $N$. Note that the right side above is continuous in $y$ and thus bounded over $K$. Here

   \begin{equation}
   \label{SigmaPW}
       \sigma(y)=
   \begin{cases}
       y & \Re y>c\\
       y{\rm{e}}^{-i\pi/2} & y\in\{\Im y> c\}\cap\{\Re y\leq c\}\\\
       y{\rm{e}}^{i\pi/2} & y\in\{\Im y< -c\}\cap\{\Re y\leq c\}
   \end{cases}
   \end{equation}

   The constants $0<c_1<1$ and $c_0>0$ were defined in \eqref{c1Const} and \eqref{estdenom}.  We note that for large $\ell$, we have $\ell!/|(y)_{\ell}|=\ell! \Gamma(y)/\Gamma(\ell+y)\sim \Gamma(y)\ell^{-y}$.  


 The bounds \eqref{FullRemEst} together with the fact that $0<c_1<1 $ show that $\mathcal{R}_{n,N,\ell}\left(\beta y,\beta\right)$ converges uniformly to zero on $K$ and Theorem \ref{choosebeta} follows.  (Of course the bound \eqref{FullRemEst}  deteriorates as the cut $\RR_-$ is approached.) 

    \end{proof} 
    
    \begin{Note}
    \label{Note12}
    {\rm To optimize the number of terms needed to achieve a desired precision of the approximation, we note that $x$ is bounded below in the complement of $\mathcal{S}_c$, and for large enough $k$ and $\ell$ not too large, the term $(2^k x/\beta)_{\ell}$ behaves like $2^{k\ell}(x/\beta)^{\ell}$ so a relatively small $\ell$ suffices to achieve high precision.}

\end{Note}

\section{Dyadic factorial expansions of various special functions}\label{Eidyadic}

\subsection{Dyadic factorial expansions of Ei in a sector containing the Stokes line}\label{EiStokes}

 For functions with a pole in the Borel plane, rather than a branch point singularity, we could integrate by parts to obtain a logarithmic branch point then apply Theorem\,\ref{ElementDecomp}. However, since no cut is needed, the techniques used in the proof of Theorem\,\ref{ElementDecomp} become simpler. We illustrate them here, and the results that we obtain: {\em global information provided by dyadic expansions} on the exponential integral, Ei, a special function often occurring in applications; see e.g. \cite{Masina} for applications and generalizations. 
 
The exponential integral is defined as $E_1(x)=\int_{x}^\infty\tfrac{e^{-t}}{t}\, dt$ on the cut plane $\CC\setminus(-\infty,0]$. $E_1$ can be analytically continued across the cut, which is a Stokes line, and we show here how this continuation can be studied numerically. See \S\ref{Eidetails} for more details about the exponential integral function and its Stokes line.

It is convenient to move the Stokes line on $\RR_+$; for this we define

\begin{equation}
\label{LiEi}
{\rm Ei}^+(x)=\mathrm{e}^{x}\int_0^{\infty e^{i0-}}\frac{\mathrm{e}^{-px}}{1-p}dp
\end{equation}
(where $0-=\phi$ is an angle with $\phi<0,\ |\phi|$ small) for $x$ with $\arg x=0^+$ and then for other $x$ by analytic continuation on the Riemann surface of the log. Note that Ei$^+$ and E$_1$ are analytic on the same Riemann surface, see \S\ref{Eidetails} for the connection between these two incarnations of the exponential integral special function.

We choose $\beta=i\pi$ so that the poles of the dyadic series for $e^{-x}$Ei$^+(x)$ are placed on $\CC\setminus i(-\infty, 0]$, cf. Theorem\,\ref{choosebeta} and Note\,\ref{choosebeta11}.

The dyadic series will then provide the function in this cut plane, unveiling numerically the Stokes phenomenon: it is known that by analytic continuation clockwise, from  $\arg(x)=0^+$ towards smaller argument, a small exponential is collected when crossing $\RR^+$ (Stokes phenomenon); upon further analytic continuation up to the cut, where $\arg x=-\pi/2^+$, the exponential becomes oscillatory, and the oscillation is revealed by the rational function expansions \eqref{dyadicEi}.

In the opposite direction, analytic continuation counterclockwise from  $\arg(x)=0^+$ to larger argument, up to the cut when $\arg x=3\pi/2^-$, unveils an asymptotic power series behavior. It is remarkable to {\em see a branch jump revealed by rational approximations}.

\subsection{Obtaining the dyadic series for $\mathrm{Ei}$}\label{dyadicfor Ei}
Using \eqref{LiEi} we proceed as in the proof of Theorem\,\ref{ElementDecomp} for $F(p)=\tfrac1{1-p}$, only here $F$ is meromorphic, so a cut is not needed in the Borel plane. 

We use the dyadic series for the Cauchy kernel \eqref{eq:deca2} for $s=1$ and $\beta=i\pi$ and we derive its dyadic series 
of $\arg x=0+$. Note that $e^{-x}{\rm Ei}^+(x)=J(0)$ where $J(t,x)$ is given by \eqref{defJdet}. 

We note that the assumption \eqref{condbeta} holds, and the proof of Theorem\,\ref{ElementDecomp} goes through with $t=0$ (and no integration in $t$).
   We obtain that the dyadic series of ${\rm Ei}$ is convergent geometrically and it is  \eqref{dyadicJdet} with $t=0$ (hence $s=1$) and $\beta=i\pi$:   \begin{equation}
  \label{eq:Eidsum}
   \mathrm{e}^{-x}{\rm Ei}^+(x)= -\sum_{m=1}^{\infty}\frac{\Gamma(m)}{2^m}\frac{1}{(y)_m}+\sum_{k=1}^\infty\sum_{m=1}^\infty\frac{\Gamma(m) \mathrm{e}^{- i\pi/2^{k}}}{(1+\mathrm{e}^{- i\pi/2^{k}})^m}\frac{1}{(2^ky)_m} \ \ \ \ \ \ (y=-i x/\pi)\qquad 
\end{equation}

\begin{Note}
\textsl{There is a dense set of poles in  \eqref{eq:Eidsum}  along $-i\RR^+$  where the dyadic expansion breaks down. (This of course does not imply \emph{actual} singularities of $\rm{Ei}^+$.)}
\end{Note}
   
When one uses the expansion for the Lerch function \eqref{RemLer} then \eqref{eq:Eidsum} can be represented as  \eqref{EiLer}. 

For approximations we need truncated series and estimates of the remainder. Writing the series \eqref{EiLer} as a sum with remainder we have
\begin{equation}
\label{dyadicEi}
 \mathrm{e}^{-x}{\rm Ei}^+(x)=-\sum_{m=1}^n \frac{\Gamma(m)}{2^m\left(\frac x{i\pi}\right)_m} -\rho_{n,0}(x)+\sum_{k=1}^{N-1} \left[ \sum_{m=1}^\ell \frac{\Gamma(m)\,\mathrm{e}^{- i\pi/2^{k}}}{\left(1+\mathrm{e}^{- i\pi/2^{k}}\right)^{m} \,  \left(\frac {2^kx}{i\pi}\right)_m}+\rho_{\ell,k}(x)\right]+R_N(x)
 \end{equation}
 where the remainders are 
\begin{equation}
\label{remEi}
\rho_{n,0}(x)=\rho_{n,0}
\left(\frac12,\frac x{i\pi}\right),\ \ \ \ \ \rho_{\ell,k}(x)= \rho_{\ell,0
}\left(\frac{e_k}{1+e_k},\frac {2^kx}{i\pi}\right),\ \ \ \ \text{ where }e_k=e^{i\pi/2^k} 
\end{equation}
with $\rho_{n,0}$ and $\rho_{\ell,0}$  
given by Lemma\,\ref{RemLer}.  From \eqref{eq:singrem} with $\beta=\pi i$,  $t=0$ we have for each $x\in \CC\setminus -i\RR^+$ and any $c\in (0,\min\{1,\frac{1}{2}\text{dist}(x,\beta(-\infty,0])\}))$  we define the contour \( \Gamma_{c} \) by \eqref{RemContGam}, from which we obtain the integral representation \eqref{RNdex}  

\begin{equation}
\label{RNdex}
R_N(0,x;\pi i)=\frac{1}{\pi i}\int_{\Gamma_{c}}  {\rm{e}}^{-qx/\pi i}\, \rho_{N}(q/\pi i,1;\pi i 
)dq
\end{equation}
This defines an analytic function in a half-plane neighborhood of $x$.  The half-plane domain may then be analytically continued by contour deformation to the region \( \Omega_c \) defined by \eqref{DomHolRem}. With $\rho_{N}$ 
given in Corollary\,\ref{dyadic Cauchy}. 

\begin{Proposition}\label{P2} 
 
   (i) For fixed $x\in \CC\setminus -i\overline{\RR^+}$ and large $n$,  $\rho_{n,0}(x)=O(2^{-n}n^{-\Im x/\pi})$. For fixed $n$ and large $x$, $\rho_{n,0}(x)=O(x^{-n})$. 

(ii) For large $l$, 
fixed $k$ and $x\in \Omega_c$, $\rho_{\ell,k}(x)=O(  |1+e_k^{-1}|^{-\ell}\ell^{-2^k\Im x/\pi})$. For fixed $\ell$ and large $2^kx$, $\rho_{\ell,k}(x)=O((2^{k}x)^{-\ell})$.

(iii) For large $N$ and fixed $x$ such that $\Re(x/ \pi i)>c$ we have 
  $|R_N(x)|\leq \frac{1}{c_0 2^{N-1} \Re(x/\pi i})$.  For the other regions $\Re(x{\rm{e}}^{-i \pi}/\pi )> c$ and $\Re(x/\pi )> c$ similar estimates follow from \eqref{eq:singrem}.
\end{Proposition}
\begin{proof}[Proof of Proposition\,\ref{P2}]

(i) For large $n$ we have, using \eqref{remEi} and \eqref{newrhonestimate2},
$$\rho_{n,0}(x)=\rho_{n,0}
\left(\frac12,\frac x{i\pi}\right) \sim \frac1{2^n}n^{ix/\pi}\Gamma(-ix/\pi)\ \ \ \ (n\to\infty)$$
while for $n$ fixed and large $x$ we use \eqref{newrhonestimate} and
$$\rho_{n,0}(x)=\frac{n!}{\left(\frac x{i\pi}\right)_{n}}\, \int_{\Gamma_{c}}\frac{ e^{-p(x+n)}}{(1+e^{-p})^{n+1}}\, dp \sim \frac{n!}{\left(\frac x{i\pi}\right)_{n}},\ \ \ \ (x\to\infty)$$

where the last estimate follows using Watson's Lemma.

(ii) For large $\ell$ and fixed $2^kx$, using \eqref{remEi} and \eqref{newrhonestimate2} we obtain
$$\rho_{\ell,k}(x)=\rho_{\ell,0 
}\left(\frac{e_k}{1+e_k},\frac {2^kx}{i\pi}\right)\sim \frac{e_k^\ell}{(1+e_k)^\ell}\, \ell^{2^kix/\pi}\, \Gamma\left(\frac{2^kx}{i\pi}\right),\ \ \ \ (\ell\to\infty)$$
while for $\ell$ fixed and large $2^kx$ we use \eqref{newrhonestimate} and  then Watson's Lemma, we obtain
$$\rho_{\ell,k}(x)= e_k^\ell\frac{\ell !}{\left(\frac{2^kx}{i\pi}\right)_\ell }\int_{\Gamma_{c}} \frac{e^{-p(x+\ell)}}{(1+e_ke^{-p})^{\ell+1}}\, dp\sim e_k^\ell\frac{\ell !}{\left(\frac{2^kx}{i\pi}\right)_\ell } ,\ \ \ \ (2^kx\to\infty) $$
We note that the contour $\Gamma_c$ used depends upon $\Re (x/\pi i)$ as defined in \eqref{RemContGam}.  This allows us to extend the asymptotics generated by Watson's lemma to the domain $\Omega_c$. 

(iii)
This is an immediate application of \eqref{eq:singrem} to \eqref{RNdex} with $\beta=\pi i$,  $t=0$ and suitable $c>0$. 
\end{proof}

\subsubsection{Numerical remarks} 
\begin{figure}[h!]
    \centering
    \includegraphics[scale=0.4]{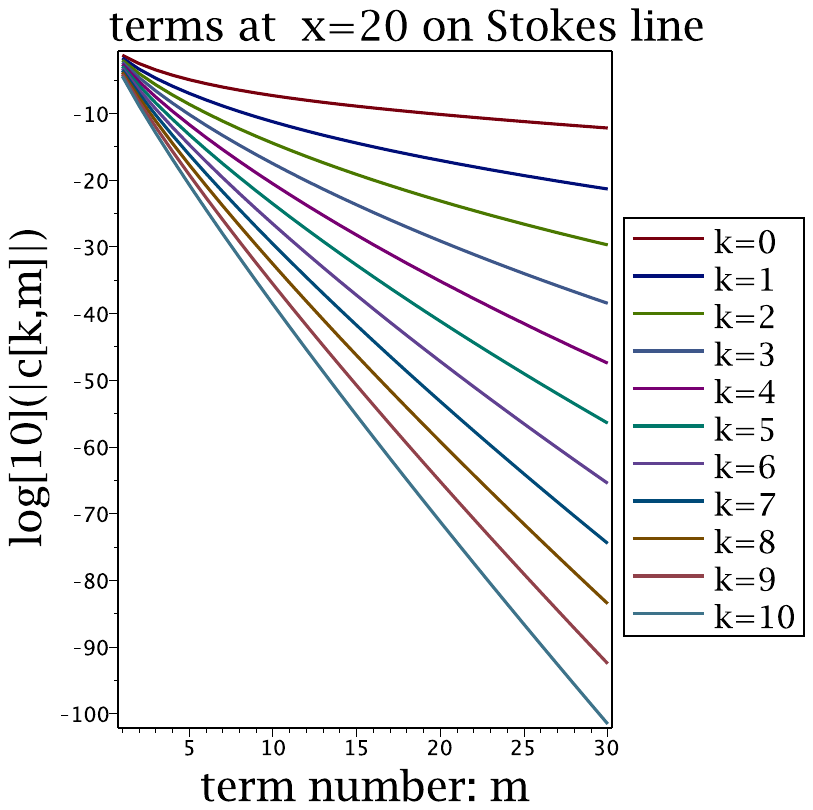}
       \caption{Size of terms in the successive series on the Stokes ray $\RR^+$ with the formula \eqref{eq:Eidsum}. This plot can be used to determine the number of terms to be kept for a given accuracy. To get $10^{-5}$ accuracy, 10 terms of the first series plus 5 from the second (with $k=1$) and so on, and all terms from the fifth series (with $k=4$) on can be discarded.}
    \label{fig12}
\end{figure}

The numerical efficiency on the Stokes line $\RR^+$, with respect to the number of terms to be kept from each of the infinitely many series in \eqref{eq:Eidsum} can be determined from  Fig. \ref{fig12}. Namely, after choosing a range of $x$ and a target accuracy, one can determine from the graphs the needed order of truncation in each individual series, as well as the number of series as described in Fig. \ref{fig12}.

In Fig. \ref{fig11} we plot the relative error in calculating Ei$^+$ on the Stokes ray.

\begin{figure}
  \centering
\includegraphics[scale=0.4]{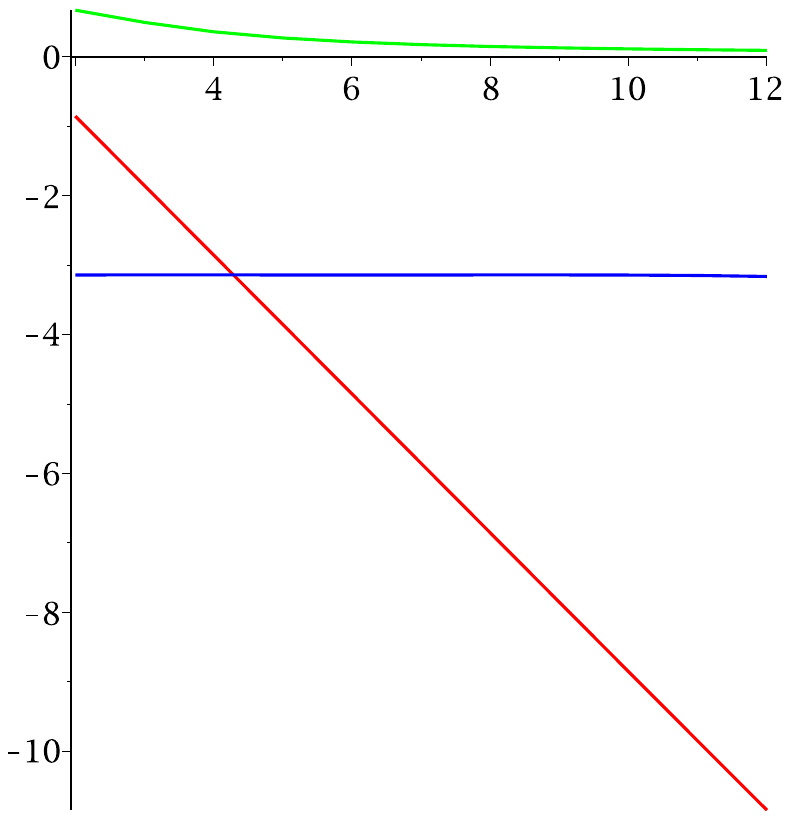}
\caption{$f(x)=\mathrm{e}^{-x}\mathrm{Ei}^+(x)$ on the Stokes line: $\Re f$, (green), $\mathrm{e}^x\Im f$, (blue), $\ln(-\Im f)$, (red), from \eqref{eq:Eidsum}. The small exponential is ``born'' on $\RR^+$, with half of the residue, as expected by comparing with $\tfrac12 \mathrm{e}^{-x}\left(\mathrm{Ei}^+(x)+\mathrm{Ei}^-(x)\right)$.  }
\label{fig15}
\end{figure}

\begin{figure}
  \centering
\includegraphics[scale=0.4]{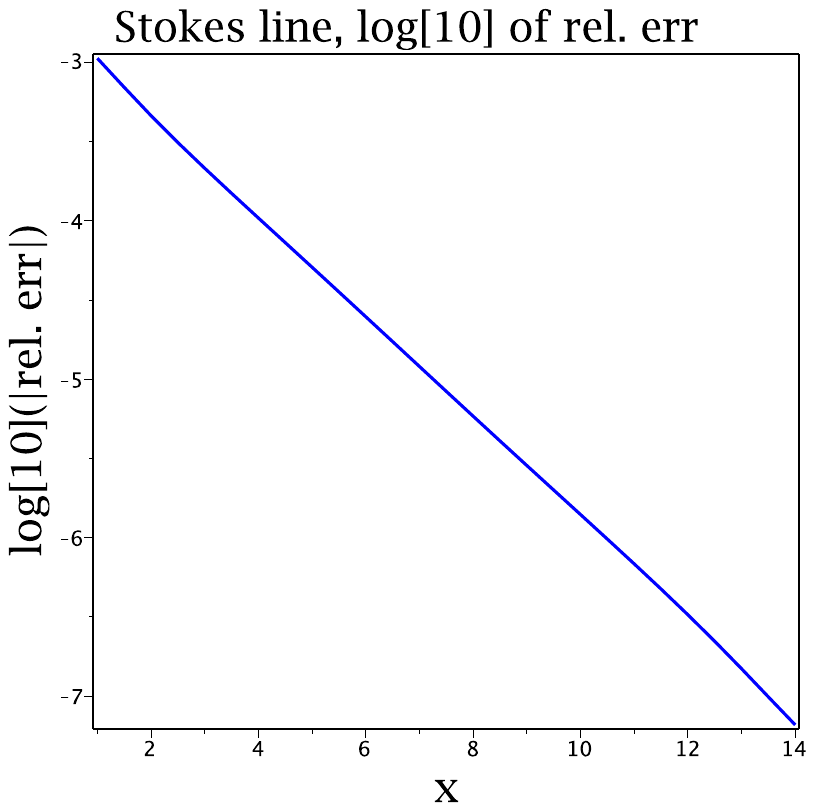}
\caption{Numerical errors for $x\in [1,14]$ for $\mathrm{e}^{-x}{\rm Ei}^+(x)$ along the Stokes line with the formula  \eqref{eq:Eidsum}. }
\label{fig11}
\end{figure}

\begin{figure}
  \centering 
\includegraphics[scale=0.3]{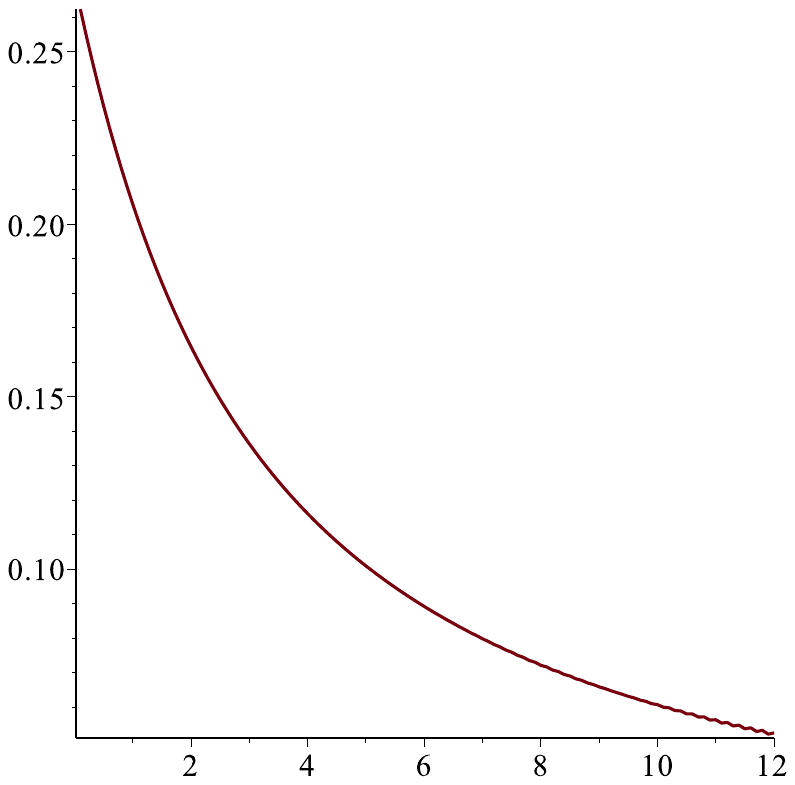} \hskip3em \includegraphics[scale=0.3]{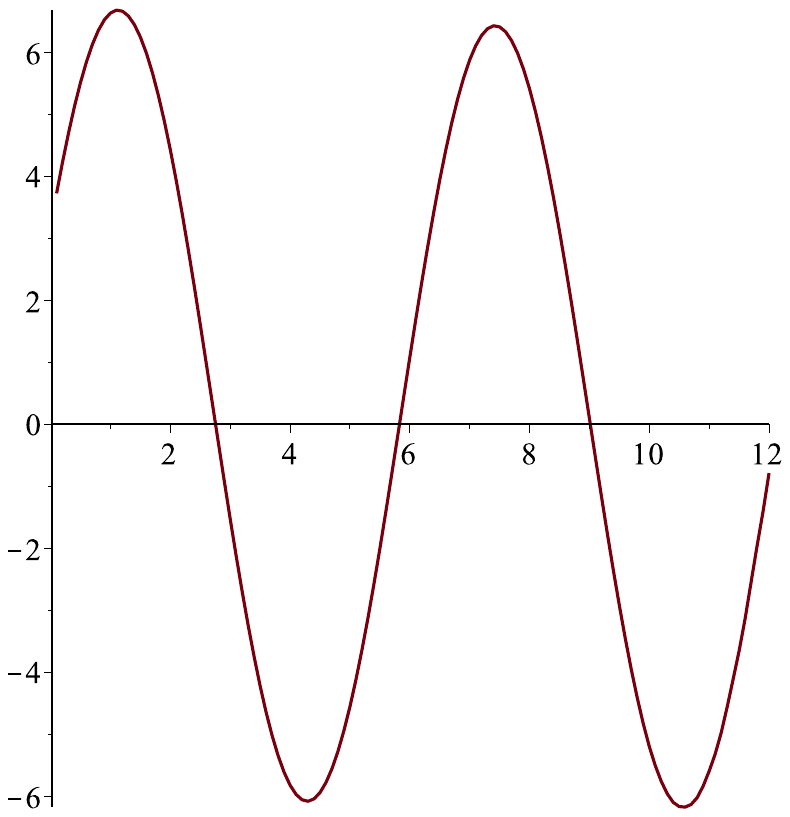} 
  \caption{The classical Stokes transition of Ei$^+$ from asymptotically decaying to oscillatory.}
\label{st1}
\end{figure}

Figure \ref{st1} below  uses the same expansion \eqref{eq:Eidsum} {for $x$} on the two sides of $-i\RR^+$; in the left picture $\Im \mathrm{e}^{-x}\mathrm{Ei}^+(x)$ is calculated for  $x\in -i\RR-0.3$ and the right one is the graph  of $\Im \mathrm{e}^{-x}\mathrm{Ei}^+(x)$  along $-i\RR+0.3$. The oscillatory behavior is due to the exponential (with amplitude $2\pi i$) collected upon  crossing the Stokes ray {$\RR^+$} (arg$\,x=-\pi/2$ is an antistokes ray for Ei$^+$).

\subsection{ Ei away from the Stokes ray, in $\CC\setminus \RR^+$}\label{Eiaway}
In \S\ref{EiStokes} we used dyadic expansions to obtain geometrically convergent expansions for Ei in $\CC\setminus i(-\infty,0]$. In this particular cut plane the convergence is least efficient due to the proximity of the antistokes line where the behavior of Ei is oscillatory. For cut planes that are away by a positive angle from the antistokes line the dyadic expansions are simpler, and more efficient.

Rotating the line of integration in \eqref{LiEi} clockwise by an angle $\pi^-$  while rotating $x$ anticlockwise by the same angle we obtain its analytic continuation as
\begin{equation}
\label{ACEip}
 {\rm e}^{-x}{\rm Ei^+}(x)=\int_0^{\infty e^{-i\pi}}\frac{e^{-px}}{1-p}\, dp=-\int_0^{\infty}\frac{e^{px}}{1+p}\, dp\ \ \ \ \ \text{ for }|\arg x-\pi|<\frac\pi 2
 \end{equation}
and we obtain its expression in terms of E$_1$ (see \S\ref{Eidetails}):
\begin{equation}
\label{E1}
{\rm e}^x{\rm Ei^+}(-x)=-\int_0^{\infty}\frac{\mathrm{e}^{-xp}}{p+1}dp=- {\rm e}^x{\rm E}_1(x) \ \ \ \ \ \text{ for }|\arg x|<\tfrac\pi 2
 \end{equation}

The dyadic expansion of \eqref{ACEip} can be obtained as in \S\ref{EiStokes}. It is convenient to take $\beta=-1$ and change variable to $y=-x$. We obtain:

\begin{Proposition}\label{PP6w} 
  The following identity holds for all $0\ne y\in\CC$ with $ |\arg\,y|<\pi/2$:
  \begin{equation}
      \label{eq:iden22}
 \mathrm{e}^y{\rm Ei}^+(-y)=\sum _{m=1}^{\infty }  \frac { (-1)^{m+1}{\mathrm{e}}
 \Gamma(m)}{ (\mathrm{e}-1)^m(y)_m} -\sum _{k=1}^{\infty }\sum _{m=1}^{\infty } \left( {\frac {\Gamma(m)\,{\mathrm{e}^{2^{-k}}}}{ \left( {
\mathrm{e}^{2^{-k}}}+1 \right) ^{m}{({2}^{k}y)_m}}}\right)
\end{equation}

Of course, in terms of Lech $\Phi$, this is the identity \eqref{eq:EiPhi}.

The remainders are similar to those of Proposition\,\ref{P2}.   
\end{Proposition}
\begin{Note}{\rm 
 The effective variable, $2^k y$, gets rapidly large for large $k$ and not many  terms of the double sum are needed in practice. Even for $y=0.1$ the first sum above (with $k=1$) requires 20 terms to give a relative error of 
 $10^{-5}$. 
    } 
\end{Note}

\subsection{Dyadic expansions for the Airy function Ai}\label{Airy}
The Airy function not only illustrates a non-trivial application of Theorem \ref{ElementDecomp}, but also shows how one can handle functions in the Borel plane which have slower decay at $\infty$.  We analyze in some detail the Airy function {\rm{Ai}} which we denote by $f(x)$ below
, as the general Bessel functions are dealt with similarly, as explained in \S\ref{genBessel}. 

After the normalization $f(x)=x^{5/4}\mathrm{e}^{-\frac23 x^{3/2}}h(4x^{3/2}/3); \ \ x=(3u/4)^{2/3}$, described in  more detail in \S\ref{sec5}, the asymptotic series of the Airy function is Borel summable:
\begin{equation}
  \label{eq:eqh1}
  h(x)=\int_0^{\infty}\mathrm{e}^{-px}F(p)dp
\end{equation}
where $F(p)={_2F_1}(1/6,5/6;1,-p)=P_{-1/6}(1+2p)$ is analytic except for a  logarithmic singularity at $-1$, see \eqref{eq:solna} and \eqref{eq:Pbranch} below. 
The decay of $F$ 
for large $p$ is relatively slow, $O(p^{-1/6})$, see \cite{nist}(15.8.2), and we integrate once by parts to improve it for Theorem \ref{ElementDecomp} to apply:
\begin{equation}
  \label{eq:intpts}
  h(x)=\frac{1}{x}+\frac1x\int_0^{\infty}\mathrm{e}^{-px}F'(p)dp
\end{equation}
where we used ${_2F_1}(1/6,5/6;1,0)=1$.  We move the singularity to $1$ by a change of variables after which we can  apply Theorem\,\ref{ElementDecomp} to \eqref{eq:Airychv} with $\beta=-1$.

\begin{equation}
  \label{eq:Airychv}
  h(-x)=-\frac{1}{x}-\frac1x\int_0^{\infty\mathrm{e}^{\pi i}}\mathrm{e}^{-px}F'(-p)dp
\end{equation}
We obtain
the dyadic series:

\begin{equation}
  \label{eq:32}
 h(-x)=-\frac{1}{x}-\sum_{m=1}^{\infty}\frac{\Gamma(m)}{ (x)_{m}}d_{m,0}-\sum_{k=1}^{\infty}\sum_{m=1}^{\infty}\frac{\Gamma(m)}{  (2^kx)_{m}}d_{m,k}
\end{equation}

Using the branch jump relation \eqref{eq:derPbranch} we see $\Delta F'(-1-t)=-iF'(t)$ and hence
\begin{equation}
  \label{dmdkm}
d_{m,0}=\frac{(-1)^{m+1}}{2\pi }\int_0^{\infty}\frac{ F'(t)\mathrm{e}^{(1+t)}}{(\mathrm{e}^{(1+t)}-1)^{m}}dt;\ \ \ \ \  d_{m,k}=-\frac1{2\pi }\int_0^{\infty}\frac{F'(t){\rm e}^{2^{-k}(1+t)}}{(\mathrm{e}^{2^{-k}(1+t)}+1)^{m}}d t
\end{equation}

Unlike in the case of Ei, the coefficients $d_m$ do not have a simple closed form expression. A convenient, and general,  way to determine them numerically is described in \S\ref{Section5}. There is an interesting expression of these coefficients in the $x$ domain: with $h$ as in \eqref{eq:eqh1} and using elementary properties of the Laplace transform we get, 
\begin{multline}
  \label{eq:sum3}
\int_0^{\infty}\frac{ F'(t)\mathrm{e}^{t+1}}{(\mathrm{e}^{t+1}-1)^{m}}dt=\sum_{j=0}^{\infty}\mathrm{e}^{-m-j+1} \binom{m+j-1}{j} \int_0^{\infty}\mathrm{e}^{-(m+j-1)t} F'(t)dt\\=\mathrm{e}^{1-m}\sum_{j=0}^{\infty}\mathrm{e}^{-j}
\binom{m+j-1}{j}[(m+j-1)h (m+j-1)-1]
\end{multline}

\begin{figure}
  \centering 
\includegraphics[scale=0.43]{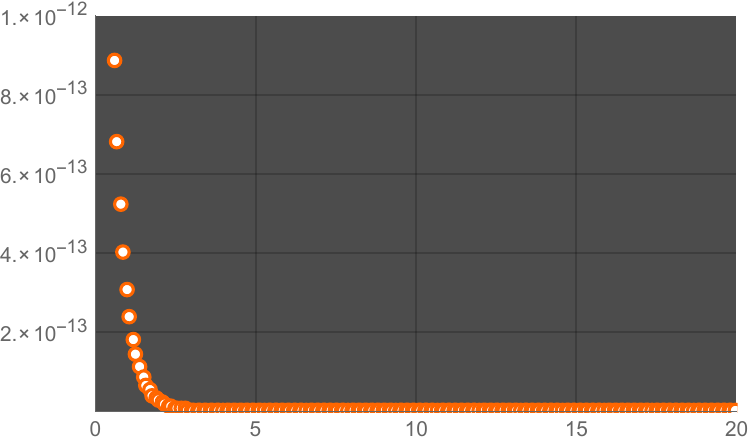}  \includegraphics[scale=0.43]{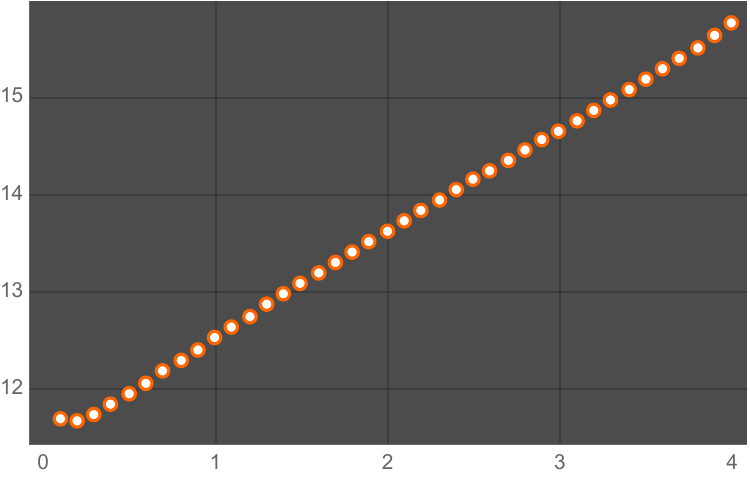} 
\caption{Relative accuracy for Ai (left), and number of exact digits (right) as functions of $x$. The total number of terms used in this calculation ranges from about 150 for small $x$  to 30 terms at $x=20$, found as explained in Fig. \ref{fig12}. The right graph plateaus at 16 digits for all  $x\ge 4$,  an artefact due to calculations being made in  Mathematica's  machine precision; thus the right graph was stopped at $x=4$. }
\label{Fig5}
\end{figure}

Fig.\,\ref{Fig5}  shows the numerical results from \eqref{eq:32} using Mathematica in machine precision to evaluate the integrals in the $d_{m,0},d_{m,k}$.

\subsection{General Bessel functions} \label{genBessel}

There are  few and relatively minor adaptations needed to deal with $K_\nu$ for more general $\nu$. After normalization, explained in \S\ref{sec5}, $F(p)$ is now the Legendre function $P_{\nu-1/2}(1+2p)$ for which the branch  jump at $-1$ is   $\Delta F(-1-p)=-2i\cos (\pi  \nu)F(p)$ (see \eqref{eq:Pbranch}) and the leading behavior at infinity is $O(p^{|\Re \nu|-1/2})$. The steps followed in the Airy case apply after integrating by parts  $k$ times until $|\Re \nu|-1/2-k<-1$.  
For $J_{\nu},Y_{\nu}$ the procedure is the same, except that the singularity is now on the imaginary line. For $J_\nu$ the singularity is on $\RR^+$ and a choice of $\beta$ as for Ei$^+$ needs to be made.

\section{Practical ways to calculate the dyadic coefficients }\label{Section5}

For a general function element $F$ the coefficients of its dyadic series are given by the integral formulas \eqref{dmk}. However, for any given $F$, these coefficients can be obtained in an efficient way as follows.

The function $F$ is represented with arbitrary accuracy by Pad\'e approximants, which we decompose by partial fractions: $F(p)\approx \sum_{i=1}^n\tfrac{c_i^{(n)}}{p-p_i^{(n)}}$. We assume, as it is typically the case, that the poles are simple. Since $F$ has only one singularity at $p=p_0$ (e.g. $p_0 =1$), the poles $p_i^{(n)}$ lie on a half line originating at $p_0$ (in our example, on $[1,+\infty)$) \cite{Stahl}. The Pad\'e {approximants}  converge in capacity; however modifications of the approximants converge uniformly at almost the same rate \cite{Stahl}. It then suffices to calculate the dyadic series for each term, $F_i^{(n)}(p)=\tfrac 1{p-p_i^{(n)}}$, which reduces to the case of the exponential integral studied in \S\ref{EiStokes}, whose dyadic coefficients are explicit. 

Due to the relatively recent, remarkable work \cite{Stahl} the accuracy of Pad\'e approximants has been established, and has an analytic expression in terms of the Green's function of a minimal capacitor, see \cite{CMP} for a summary of these fundamental results and practical algorithms to estimate this accuracy.

\section{Dyadic expansions for the Psi function and a curious identity}\label{Psi}

The digamma function, or Psi function, defined as
$$ \Psi(z)= \frac{\Gamma'(z)}{\Gamma(z)}$$
is a meromorphic function with simple poles of residue $-1$ at $z=-1,-2,-3,\ldots$ see \cite{nist}(5.2.1).

We find, and state in Proposition\,\ref{curious} the dyadic series for the Psi function, and one for differences of Psi functions, which yields a curious identity, \eqref{eq:ide}, which appears to be new.

\begin{Proposition}\label{curious}

(i) We have, for all $x\in\CC\setminus(-\infty,0]$,
\begin{equation}
  \label{eq:doublexp}
  \Psi(x+1)=\ln x-\sum_{k=1}^{\infty}\Phi(-1,1,2^kx+1) =\ln x+\sum_{k=1}^{\infty}\sum_{j=1}^{\infty}{\frac{(j-1)!}{2^j(2^k x+1)_j}}
\end{equation}

(ii)  For all $x\in\CC\setminus(-\infty,0]$ we have
\begin{equation}
  \label{eq:modif}
\frac{1}{2} \Psi \left(\frac{x}{2}+\frac{1}{2}\right)-\frac{1}{2} \Psi
   \left(\frac{x}{2}\right)= \frac{1}{2x} \ - \ \frac{1}{2^2(x)_2 }+\cdots+\frac{(-1)^{n-1}\Gamma(n)}{2^n(x)_n}+\cdots
\end{equation}

For $\Re x>0$ we have
\begin{equation}
  \label{eq:modif2}
  \frac{1}{2} \Psi \left(\frac{x}{2}+\frac{1}{2}\right)-\frac{1}{2} \Psi
   \left(\frac{x}{2}\right)= \int_0^1\frac{t^{x-1}}{t+1}dt
   \end{equation}

(iii) The identity \eqref{eq:ide} holds.

\end {Proposition}

\begin{proof}

(i) Replacing $p$ by $-p$ in \eqref{eq:deca1} we get
\begin{equation}
  \frac1p-\frac{1}{\mathrm{e}^p-1}=\sum_{k=1}^\infty\frac {{\rm e}^{-\frac {p}{{2}^{k}}}}{2^k\left( {{\rm e}^{
-{\frac {p}{{2}^{k}}}}}+1 \right)} 
\end{equation}
On the other hand  we have, see \cite{Book} eq. (4.61) p. 99,
\begin{equation}
  \label{eq:lng}
  \frac{\Gamma'(x+1)}{\Gamma(x+1)}-\ln x=\int_0^\infty\left(\frac1p-\frac{1}{\mathrm{e}^p-1}\right)\mathrm{e}^{-xp}dp
\end{equation}
Thus, changing the variable of integration to $q=p/2^k$ we get
\begin{equation}
  \label{eq:gamma2}
 \Psi(x+1)= \frac{\Gamma'(x+1)}{\Gamma(x+1)}=\ln x+\sum_{k=1}^{\infty}\int_0^\infty\frac{\mathrm{e}^{-q(1+2^k x)}}{1+\mathrm{e}^{-q}}\, dq
\end{equation}
Using \eqref{ClassicalLer} together with \eqref{factorialseriesforPhi} and taking $z=1/2$, 
 we obtain the dyadic factorial expansion \eqref{eq:doublexp}.

(ii) Consider the functional equation
\begin{equation}
  \label{eq:difdif1}
  f(x+1)+f(x)=\frac 1x  
\end{equation}
After Borel transform (i.e. substituting \eqref{fisLF} in \eqref{eq:difdif1}) we obtain $(\mathrm{e}^{-p}+1)F(p)=1$, yielding
\begin{equation}
  \label{eq:id4}
  f(x)=\int_0^\infty\,\frac{\mathrm{e}^{-px}}{\mathrm{e}^{-p}+1}\, dp= \int_0^\infty\,\sum_{n=0}^\infty\,(-1)^n\,\mathrm{e}^{-p(x+n)}\, dp = \sum_{n=0}^\infty \frac{(-1)^n}{x+n}
\end{equation}
where the interchange of summation and integration is justified, say, by the monotone convergence theorem applied to $\sum_{n=0}^{2N}(-1)^n\mathrm{e}^{-p(x+n)}$. Of course, the integral converges only for $\Re x>0$, but the series converges for all $x\notin\{0,-1,-2,...\}$. Therefore $f(x)$ is meromorphic, having simple poles at $x=-n,\ n\in\NN$.

 On the other hand $f(x)=\tfrac12\psi(\tfrac x2+\tfrac12)-\tfrac12\psi(\tfrac x2)$ which follows from integrating the  identity
$$\psi'(z)= \sum_{n=0}^{\infty}(z+n)^{-2}$$
 (see \cite{Ahlfors},  (31) p. 200) between $z=\frac x2$ and $z=\frac{x+1}2$.
 
 The integral representation \eqref{eq:modif2}
 then follows by substituting $\mathrm{e}^{-p}=t$ in \eqref{eq:id4} and the factorial expansion in \eqref{eq:modif} is then obtained as usual, by integration by parts.

 (iii) is obtained straightforwardly by combining \eqref{eq:modif} with \eqref{eq:doublexp}.
 
 \end{proof}

 \section{Duplication formulas and incomplete Gamma functions}\label{dupl}
 Some applications, such as the ones in \S\ref{erfc} and \S\ref{op}, require fractional powers.
 In this section Lemma\,\ref{L1L} generalizes Lemma \ref{L1} to fractional powers of $p$ and find simple dyadic representations for some other classes of special functions. We find factorial series with coefficients having closed form expressions in terms of polylogarithms.

Recall that the polylog is defined as
$$  \mathrm{Li}_{s}\left(z\right)=\sum_{k=1}^\infty \frac{z^k}{k^s}$$
for any $s\in\CC$. The series converges for $|z|<1$ and $\mathrm{Li}_{s}(z)$ is defined by analytic continuation for other values of $z$, \cite{nist}25.12.10. It has the integral representation	
\begin{equation}
  \label{eq:intrep}
  \mathrm{Li}_{s}\left(z\right)=\frac{z}{\Gamma
\left(s\right)}\int_{0}^{\infty}\frac{x^{s-1}}{\mathrm{e}^{x}-z}dx
\end{equation}
when $\Re s>0$ and $\arg(1-z)<\pi$, or $\Re s>1$ and $z=1$, \cite{nist}25.12.11.

$ \mathrm{Li}_{s}(z)$ satisfies the general duplication formula
\begin{equation}
  \label{eq:dupl}
  f(z)+f(-z)=2^{1-s}f(z^2)
\end{equation}
(see \cite{Lewin}; also, \eqref{eq:intrep}, \eqref{eq:dupl} are easily checked directly).

\begin{Lemma}[A ramified generalization of \eqref{eq:deca1}]\label{L1L}{\rm 
 \textsl{ The following identity holds in $\CC\setminus \{0\}$:
\begin{equation}
  \label{eq:li2}
 \pi p^{s-1}=\Gamma(s)\sin(\pi s)\left[ \mathrm{Li}_s\left(\mathrm{e}^{-p}\right) - \sum _{k=1}^{\infty} 2^{-k (1-s)} \mathrm{Li}_s\left(-\mathrm{e}^{-2^{-k} p}\right)\right] \ \ \ \ \text{ if } s<1
\end{equation}
which reduces to \eqref{eq:deca1} if $s=0$.
}}
\end{Lemma}

\begin{proof}
 Let $s<1$. As in the proof of Lemma\,\ref{L1} we iterate \eqref{eq:dupl} $n$ times:
 \begin{multline}
 f(z)=-f(-z)+2^{1-s}f(z^2)=-f(-z)-2^{1-s}f(-z^2)+2^{2(1-s)}f(z^4)=\ldots \\
 =2^{n(1-s)}f\left(z^{2^n}\right)-\sum_{j=0}^{n-1}  2^{j(1-s)} f\left(-z^{2^j}\right)
  \end{multline}
 where, taking $z=e^{p/2^n}$ we obtain

\begin{equation}
  \label{eq:li222}
 2^{-n
   (1-s)}\text{Li}_s\left(\mathrm{e}^{-p/2^n}\right)= \text{Li}_s\left(\mathrm{e}^{-z}\right) - \sum _{k=1}^n 2^{-k (1-s)} \text{Li}_s\left(-\mathrm{e}^{-p/2^k}\right)
\end{equation}
 We use the following identity \cite{nist}(25.12.12):
$$\mathop{\mathrm{Li}_{s}\/}\nolimits\!\left(u\right)=\mathop{\Gamma\/}\nolimits%
\!\left(1-s\right)\left(\mathop{\ln\/}\nolimits\frac{1}{u}\right)^{s-1}+\sum_{%
j=0}^{\infty}\mathop{\zeta\/}\nolimits\!\left(s-j\right)\frac{(\mathop{\ln\/}%
\nolimits u)^{j}}{j!},\ \ \ \ \ \ \ \ \ \ \ \ s\ne 1,2,\ldots,\ \ \ \ |\ln u|<2\pi$$
which, for $u=\mathrm{e}^{-p/2^n}$ becomes
$$2^{-n(1-s)}\text{Li}_s\left(\mathrm{e}^{-p/2^n}\right)=\Gamma(1-s) p^{s-1}+  2^{-n(1-s)}\sum_{j=0}^\infty  \mathop{\zeta\/}\nolimits\!\left(s-j\right)\frac{(\mathop{\ln\/}%
\nolimits u)^{j}}{j!},$$
where we see that, for $s<1$, $\lim_{n\to\infty}2^{-n(1-s)}\text{Li}_s\left(\mathrm{e}^{-p/2^n}\right)=\Gamma(1-s) p^{s-1}$.
Thus,  taking the limit $n\to\infty$ in \eqref{eq:li222} we obtain
$$ z^{s-1}\Gamma(1-s)= \text{Li}_s\left(\mathrm{e}^{-z}\right) - \sum _{k=1}^{\infty} 2^{-k (1-s)} \text{Li}_s\left(-\mathrm{e}^{-2^{-k} z}\right)  $$
from which \eqref{eq:li2} follows by using the reflection formula $\Gamma(s)\Gamma(1-s)=\pi/\sin(\pi s)$. 
\end{proof}

\subsection{Dyadic series for incomplete gamma functions and erfc}\label{erfc} 
The incomplete gamma function, which arises as solution to various mathematical problems, is defined by
$$\Gamma(s,x)=\int_x^\infty \, t^{s-1}\, \mathrm{e}^{-t}\, dt$$
and has as a special case the error function, 
$$\text{erfc}(x)=\frac{2}{\sqrt{\pi}}\int_x^\infty\, \mathrm{e}^{-t^2}\, dt=\frac{1}{\sqrt{\pi}}\, \Gamma\left(\frac 12,x^2\right)$$ 
Noting that
$$\int_0^{\infty}(1+p)^{s-1}\mathrm{e}^{-xp}dp=\mathrm{e}^x x^{-s}\Gamma(s,x)$$
we see that $\mathrm{e}^x x^{-s}\Gamma(s,x)$ is the Laplace transform of a function which has a ramified singularity if $s\not\in\ZZ$. In this case we apply Lemma\,\ref{L1L} and obtain the expansion, for $s<1$
\begin{equation}
  \label{eq:gammat}
  \Gamma(1-s)\mathrm{e}^xx^{-s}\Gamma(s,x)=\mathcal{L} \,\text{Li}_s\left(\mathrm{e}^{-p-1}\right)-\sum _{k=1}^{\infty} 2^{-k (1-s)} \mathcal{L}\,\text{Li}_s\left(-\mathrm{e}^{-2^{-k} (p+1)}\right)
\end{equation}
and in particular
\begin{equation}
  \label{eq:erfc}
   \pi e^x x^{-1/2}\text{erfc}\left(\sqrt{x}\right)=\mathcal{L}\,\text{Li}_{\frac{1}{2}}\left(\mathrm{e}^{-p-1}\right)-\sum
   _{k=1}^{\infty } 2^{-k/2} \mathcal{L}\,\text{Li}_{\frac{1}{2}}\left(-\mathrm{e}^{-2^{-k}
   (p+1)}   \right)
\end{equation}
From this point on, the dyadic expansions are obtained by calculating the factorial expansion of each term in \eqref{eq:erfc}. For example, the first term in \eqref{eq:erfc} has the factorial series
\begin{equation}
\label{Lli}
 \mathcal{L}\,\text{Li}_{\frac{1}{2}}\left(\mathrm{e}^{-p-1}\right)=\int_0^1t^{x-1}\text{Li}_{\frac{1}{2}}\left(\frac t{\mathrm{e}}\right)dt=\sum_{k=0}^\infty\frac{(-1)^k}{\mathrm{e}^k(x)_{k+1}}\text{Li}_{\frac{1}{2}}^{(k)}\left(\mathrm{e}^{-1}\right):=\sum_{k=0}^\infty\frac{c_k}{(x)_{k+1}} 
 \end{equation}
with
\begin{equation}
\label{Llicoeffs} 
c_k =(-1)^k\sum_{j=0}^k s(k,j)\text{Li}_{\frac{1}{2}-j}\left(\mathrm{e}^{-1}\right)
 \end{equation}
where $s(k,j)$ are the Stirling numbers of the first kind (see \S\ref{form} for details), where we used the formula
\begin{equation}\label{derLi}
\frac{d^k}{dz^k}\text{Li}_{\nu}\left(z \right)=z^{-k}\sum_{j=0}^k s(k,j)\text{Li}_{\nu-j}(z)
\end{equation}
We note that \eqref{derLi} can be verified by an inductive argument on the relation $\frac{d}{dz}\text{Li}_s(z)=z^{-1}\text{Li}_{s-1}(z)$ and we used the formulas for $s(k,j)$ discussed in \S\ref{form}.

\section{Dyadic resolvent identities}\label{op}
Dyadic decompositions translate into representations of the resolvent of a self-adjoint operator in a series involving the unitary evolution operator at specific discrete times:
\begin{Proposition}\label{C8}\textsl{ 
(i) Let $\mathcal{H}$ be a Hilbert space,  and $A$ a bounded or unbounded self-adjoint operator. Let $U$ be the unitary evolution operator generated by $A$, $U_t=\mathrm{e}^{-itA}$. If $\lambda\in\RR^+$, then}
\begin{equation}
  \label{eq:resol1}
(A-i\lambda)^{-1}  =i(1-\mathrm{e}^{-\lambda}U_1)^{-1}-i\sum_{k=1}^{\infty}\frac1{2^k}(1+e^{-{\lambda}/2^{k}}U_{2^{-k}})^{-1} 
\end{equation}
and \eqref{eq:resol2} follows.

\textsl{{Convergence holds in the strong operator topology.}  For $\lambda<0$ one simply complex conjugates \eqref{eq:resol1}.  (The limits cannot, generally, be interchanged.)}

 (ii) Assume $A$ is a  positive operator (thus self-adjoint) and  $0\notin\sigma(A)$. Let $T_t$ be the semigroup generated by $A$, $T_t=\mathrm{e}^{-tA}$. Then 
 \begin{equation}
   \label{eq:Semigr}
A^{-1}=(1-T_1)^{-1}-\sum_{k=1}^{\infty}2^{-k}(1+T_{1/2^k})^{-1}=\sum_{j=1}^{\infty}T_{j}-\lim_{\ell\to\infty}\sum_{k=1}^{\ell}\sum_{j=1}^{\infty}2^{-k}(-1)^jT_{j/2^k}
 \end{equation}
 
 \textsl{where now convergence is in  operator norm. More generally, for $ s<1$, $s\notin\ZZ$,}
\begin{equation}
  \label{eq:li21}
 \pi A^{s-1}=\Gamma(s)\sin(\pi s)\left[ \mathrm{Li}_s\left(T_1\right) - \sum _{k=1}^{\infty} 2^{-k (1-s)} \mathrm{Li}_s\left(-T_{1/2^k}\right)\right]
\end{equation}
\textsl{in operator norm} 
\end{Proposition}
For a discussion of the polylog function $\mathrm{Li}_s(z)$ see \S \ref{dupl}.

\begin{proof}
(i) We recall the projection-valued measure  spectral theorem for  self-adjoint operators. If $\mathcal{H}$ and $A$ are as above  and  $g:\RR\to\RR$ is a  Borel function (or  a complex one, by writing $g=g_1+ig_2$), then $g(A)=\int_{-\infty}^{\infty}g(q) dP_{q}$ where $\{P_{\Omega}\}$ are the projection-valued measures induced by $A$ on $\sigma(A)$ 
 (see \cite{ReedAndSimon} Theorem VIII.6 p. 263). The spectral theorem together with  \eqref{eq:decx3}   for  $p=\lambda+iq$ give 
 \begin{equation}
   \label{eq:decop}
   (1-\mathrm{e}^{-\lambda}U_1)^{-1}-\sum_{k=1}^n2^{-k}(1+ \mathrm{e}^{-2^{-k}\lambda}U_{2^{-k}})^{-1}=\epsilon(1-\mathrm{e}^{-\lambda\epsilon}\mathrm{e}^{-i\epsilon A})^{-1}=\int_{\RR}\frac{\epsilon dP_q}{1-\mathrm{e}^{-\epsilon(\lambda+iq)}}
 \end{equation}
where  $\epsilon_n:=\epsilon=2^{-n}$. An elementary calculation shows that the modulus of the integrand is uniformly bounded by $\lambda^{-1}$. Since the integrand converges pointwise to $(\lambda+iq)^{-1}$ as $\epsilon\to 0$, dominated convergence shows that the integral converges to $(\lambda+iA)^{-1}$. Dominated convergence  also shows that the integrand, seen as a multiplication  operator, converges in the strong operator topology, implying the result.

(ii) The proof,  based on the same argument as in (i), is simpler  and we omit it. For \eqref{eq:li21} we combine this argument with  Lemma \ref{L1L} below. The sums in \eqref{eq:Semigr} are manifestly convergent in the operator norm since  $\|T_t\|<1$ and $T_t>0$.
  \end{proof}

\section{Dyadic series of typical functions occurring in applications; resurgence}\label{Resfun} 

Generic systems of meromorphic ODEs, difference equations and other classes of problems commonly occurring in applications have solutions characterized by a special Borel plane structure. Their Borel transforms\footnote{Recall that if $f(x)=\int_0^\infty e^{-px}F(p)dp$ we say that $f$ is the Laplace transform of $F$, and $F$ is the inverse Laplace trasform, or Borel transform, of $f$.} satisfy the following conditions:

 \ \ \ \ {\em(A1)}\label{A1} they have at most exponential growth at infinity (meaning a finite exponentially weighted $L^1$ norm, see \eqref{nubounded} and \S\ref{PT}(ii) for a precise formulation

  \ \ \ \ {\em (A2)} their singularities are equally spaced \footnote{Equal spacing of singularities is a typical characteristic of resurgent functions. As keenly pointed out by one of the referees, the results in \S \ref{Resfun} can be easily extended to the case where the singular set is of the form $\{n_j\lambda_k\}_{k=1, j \in \mathbb{N}}^n$, where $\{n_j\}_j$ is an  increasing sequence of positive numbers such that $\sum_j \rho_k^{n_j}<\infty$, where $\rho_k={\rm{e}}^{-(\mu-\nu)c_k |\lambda_k|}<1$ and $\lambda_k$ are $\mathbb{Z}$-independent complex numbers, as discussed in \S \ref{PT}.    } along finitely many rays, and

  \ \ \ \ {\em (A3)} at each singularity there exists locally Frobenius-type convergent expansions in fractional powers and possibly logs (see e.g. \cite{Duke,Braaksma});  the singularities are integrable.\footnote{The $L^1_{\rm loc}$ nature of the power can be often arranged by a suitable substitution. Alternatively, integrals through the $\omega_i$'s can be replaced by integrals avoiding the singularities, see \S\ref{PT}, \eqref{defR}.}

   \ \ \ \ {\em (A4)}\label{A4} the singularities are non-resonant (this is the generic case, see \S\ref{PT}(i) for a precise formulation)
 
In fact, more is true for the aforementioned  solutions: the singularities on the Riemann surface of solutions of the same equation are interconnected in an explicit fashion, and possess a set of deep properties --they are resurgent in the sense of \'Ecalle, see \cite{Ecalle,Duke}. 

\begin{Definition}
We say that a function $F(p)$ satisfies  {\em Assumption {\bf (A)} }  if it has the properties {(A1)}-{(A4)}.
\end{Definition}

\subsection{Decomposition of resurgent functions into function elements}
 We defined {\em function elements} to be
resurgent functions with only one regular singularity on the first Riemann sheet, and with algebraic decay at infinity, see \S\ref{FuncEl}.
There are two main properties of function elements which do not  hold for general resurgent functions: decay at infinity and the property of having only one singularity. Resurgent functions can be nonetheless decomposed into {\emph{function elements}}. 

 To avoid cumbersome details and keep the presentation clear, we present the essential steps and formulate  Theorem\,\ref{TT1} for the case where the resurgent function has the form encountered as solutions of generic meromorphic ODEs.

\begin{Theorem}\label{TT1}
 The Laplace transform $f$ of functions $F$ satisfying Assumption {\bf (A)} can be written, after a  translation of the variable, as a geometrically convergent series of Laplace transforms of function elements plus an entire function. 
\end{Theorem}

\begin{Note}  The exponential integral and the $\Psi$ function treated in \S\ref{Psi} are examples of elements with nonramified singularities. Airy and Bessel functions treated in \S\ref{Airy} and \S\ref{genBessel} are examples of elements with ramified singularities, treated via the Cauchy kernel decomposition. The incomplete gamma function and the error function  treated in \S\ref{erfc} have power-ramified singularities for which a polylog dyadic expansion (Lemma\,\ref{L1L}) gives more explicit decompositions. Theorem\,\ref{TT1} extends these techniques to general resurgent functions.
\end{Note}

\subsection{Proof of Theorem \ref{TT1}}

The decomposition of $F$ is constructed in \S\ref{PT}, then in \S\ref{pfTT1} it is proved that this decomposition has the desired properties.

\subsubsection{Decomposition in function elements} \label{PT}
In this section we describe how a resurgent function can be decomposed into function elements.

Let $F(p)$ satisfy Assumption {\bf(A)}. We introduce the following notations.

\ \ \ \ \ $(i)$\ \  \ Denoting  by $\omega_i$ the singularities of $F$, then by {\em(A2)}, each $\omega_i$ is of the form $j\lambda_k$, with $j\in\ZZ^+$ and $\lambda_k\in\{\lambda_1,\ldots,\lambda_n\} $.

By {\em(A4)} they are assumed non-resonant, in the sense that
$\lambda_1,\ldots,\lambda_n$ are linearly independent over $\ZZ$ and of different complex arguments;

\ \ \ \ \ $(ii)$\ \  We define the space $\mathscr{P}$ of smooth curves starting at the origin, traveling forward towards infinity while avoiding the singularities: if $\Omega\coloneqq\CC\setminus \bigcup_{j=1}^n\lambda_j \NN$, let
\begin{equation}
    \mathscr{P} \coloneqq\left\{\gamma:(0,1)\to \Omega \; \bigg| \; \gamma(0^+)=0, \;\gamma \text{ smooth }, \frac{d}{dt}|\gamma(t)|>0, \; \arg\left(\gamma(t)\right)\; \text{is monotone}, \; \lim_{t\to 1^-}|\gamma(t)|=\infty \right\}
\end{equation}
For every $\gamma\in \mathscr{P}$ we adjoin to $\mathscr{P}$ the path $-\gamma$ with reversed orientation; traveling from infinity to the origin rather than to infinity from the origin.  We denote the larger set of curves by $\mathscr{P^{\pm}}$.

We further restrict our attention and consider only the paths whose length does not grow too fast: for $\kappa>0$ let
\begin{equation}\label{defR}
\mathscr{R} \equiv   \mathscr{R}_\kappa \coloneqq\left\{\gamma \in \mathscr{P}^\pm: \text{so that for all} \;  R>0, \ \ len\left(\gamma\cap\DD_R\right) < \kappa R\right\}
\end{equation}
where $len\left(\gamma\cap\DD_R\right)$ denotes the length of the part of $\gamma$ contained in the disk $\DD_R$.

\ \ \ \ \ $(iii)$ \ \ By assumption {\em(A1)} there is a $\nu>0$ such that the following sup is bounded:

\begin{equation}\label{nubounded}
\|F\|_{1,\nu} : =\sup_{\gamma\in \mathscr{R}}\int_{\gamma}|F(p)|\mathrm{e}^{-\nu |p|}|dp|
<\infty
\end{equation}
Sometimes it will be convenient to compute a weighted $L^1$ norm along a specific ray with direction $\varphi$.  For this we use the notation for the function space and norm respectively:

\begin{equation}\label{weightedL1direcphi}
   L^1_{\nu,\varphi}:=L^1\left({\rm{e}}^{i\phi}\RR_+, {\rm{e}}^{-\nu|s|}|ds| \right) , \qquad \|F\|_{1,\nu,\varphi}:=\int_{0}^\infty |F(x{\rm{e}}^{i\varphi})|{\rm{e}}^{-\nu x} dx 
\end{equation}

\ \ \ \ \ $(iv)$ We denote by $S_i$ thin, non-intersecting half-strips containing exactly one singularity $\omega_i$. More precisely, in the notation {\em(i)}, each $\omega_i\equiv\omega_{jk}=j\lambda_k$, and $S_i\equiv S_{jk}$. It is these half-strips that determine our choice of branch cut; we take the cuts to be along the axis of symmetry of each $S_{jk}$. We can assume, for simplicity, that all $S_{jk}$ with the same $k$ are translates of each other and that each strip is eventually bordered by straight lines of arguments $\theta_k$. For good convergence (as in Lemma\,\ref{GeomLemmaRes}) the directions of $S_{jk}$ cannot be orthogonal to $\lambda_k$, and of course, they cannot be parallel either: we assume 
\begin{equation}\label{assumethetak}
0<\delta_{1,k}<|\theta_k-\arg\lambda_k|<\delta_{2,k}<\frac\pi2
\end{equation}
We will assume $\delta_{2,k}$ small enough.

\begin{figure}
\centering
\hspace{2cm}\includegraphics[scale=0.25]{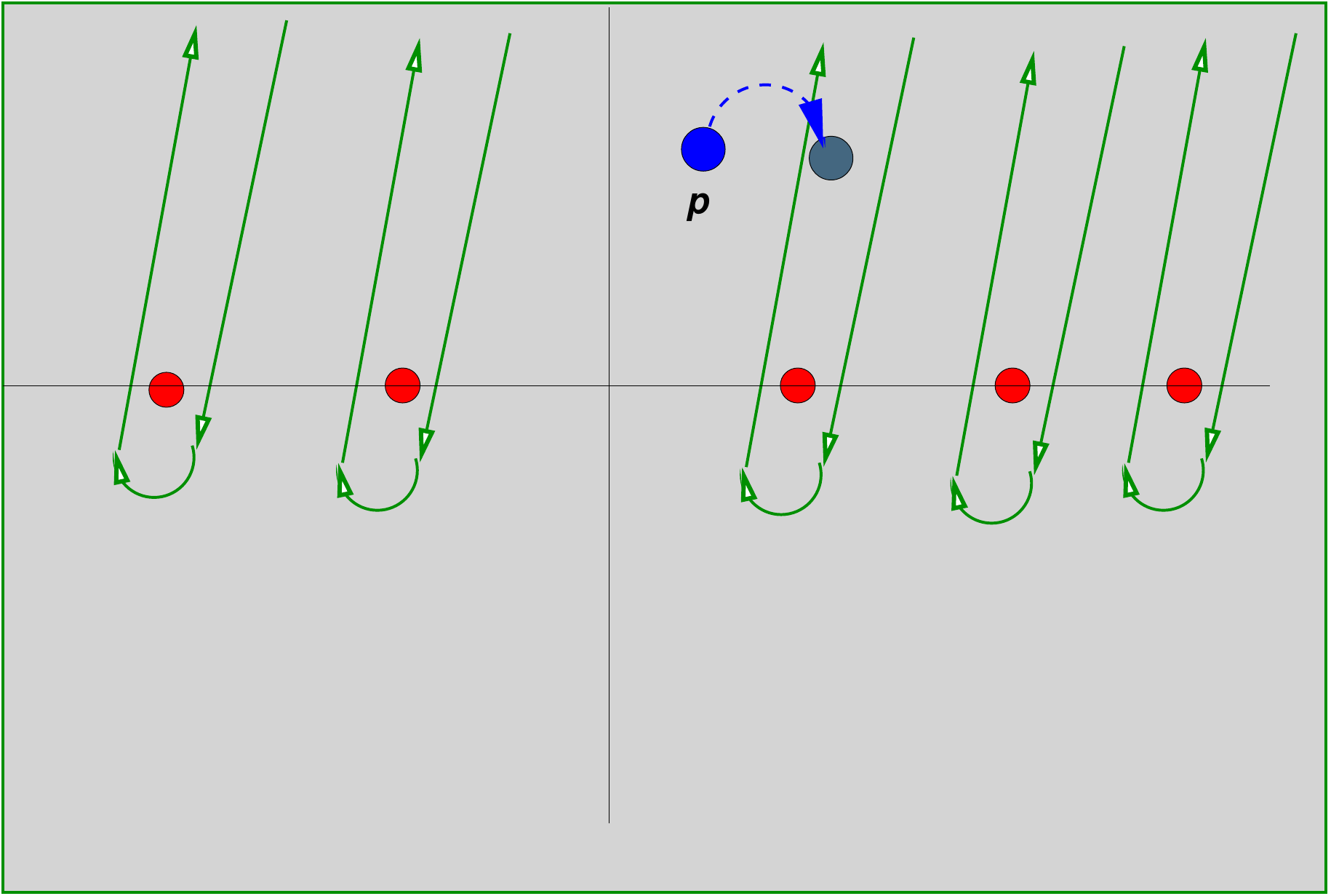}
\vskip -0.5cm
\caption{ The contours $\mathcal{C}_i$.}
\label{fig16}
\end{figure}

 We let $\mathcal{C}_i=\partial S_i$, see Fig. \ref{fig16}, non-intersecting H\"ankel contours around the $\omega_i$, going 

 towards $\infty$; $\mathcal{C}_i$ are traversed anticlockwise. 
 
  Let $\mathcal{A}$ be the complement of the union of $S_i$.

Let
\begin{equation}
\label{eq:decp}
G(p)=F(p)-\sum_{\omega_i} \frac{\exp(\mu_i p)}{2\pi i}\int\limits_{\mathcal{C}_i} \frac{F(s)\exp(-\mu_i s)}{s-p}ds
\end{equation}

 Define 
 \begin{equation}
  \label{eq:F_i}
  F_i(p)=\frac{\exp(\mu_i p)}{2\pi i}\int\limits_{\mathcal{C}_i} \frac{F(s)\exp(-\mu_i s)}{s-p}ds\, \equiv\, {\exp(\mu_i p)}\, \tilde{F}_i(p)
\end{equation}
and let
 \begin{equation}
    \label{eq:decp}
    G(p)=F(p)-\sum_{\omega_i} F_i(p)
      \end{equation}
where:

\ \ \ \ \ $(iv_1)$\ \   $|\mu_i|=\mu>\nu$, 

\ \ \ \ $(iv_2)$\ \   {If $\omega_i=j\lambda_k\equiv\omega_{jk}$, we choose $\arg(\mu_i)$ so that $\Re(\mu_i s)>0$ for large $s$ on $\mathcal{C}_i$; this means that $\mu_i\equiv\mu_{jk}=\mu e^{-i\theta_k}$ with $\theta_k$ satisfying \eqref{assumethetak}.  We will sometimes sum over the index $j\in\NN$ only. Note that the function defined in \eqref{eq:F_i} is $F_i\equiv F_{jk}$.

\subsubsection{Proof of Theorem \ref{TT1}}\label{pfTT1}
\begin{proof}
The main steps in the proof are as follows. 
We show that each $\tilde{F}_i$ defined in \eqref{eq:F_i} is a function element: Lemma\,\ref{singFi} shows that it has only one singularity on the first Riemann sheet and Lemma\,\ref{decaytFi} proves algebraic decay at infinity.
Lemma\,\ref{GeomLemmaRes} proves the series in \eqref{eq:decp} converges uniformly on compact sets, hence $G$ is analytic. Finally, Lemmas\,\ref{FiL1nu}, and \ref{GisLaplace} show that the entire function $G$ defined by \eqref{eq:decp} satisfies necessary growth conditions (it is at most exponentially increasing, so that it is Laplace transformable) for Theorem \ref{TT1} to hold. 

\

\begin{Lemma}\label{singFi}
  On the first Riemann sheet, each $F_i$ in \eqref{eq:F_i} has precisely one singularity, namely at $\omega_i$. Furthermore $F-F_i$ is analytic at $\omega_i$.
\end{Lemma}
\begin{proof}
  Let  $p\ne \omega_i$.  If $p$ is outside $\mathcal{C}_i$ then  function $F_i$ is manifestly analytic at $p$. To analytically continue in $p$ to the interior of $\mathcal{C}_i$  it is convenient to first deform $\mathcal{C}_i$ past $p$, collecting the residue. We get
  \begin{multline}\label{eq:ResidueResurg}
      F_i(p)=\frac{\exp(\mu_i p)}{2\pi i}\left[ \int\limits_{\widetilde{\mathcal{C}}_i} \frac{F(s)\exp(-\mu_i s)}{s-p}ds+2\pi i F(p) \exp(-\mu_i p)\right]\\
      =F(p)+\frac{\exp(\mu_i p)}{2\pi i}\int\limits_{\widetilde{\mathcal{C}}_i} \frac{F(s)\exp(-\mu_i s)}{s-p}ds
  \end{multline}
where now $p$ sits inside $\widetilde{\mathcal{C}}_i$, and the new integral is again manifestly analytic. 

Thus $F_i$ is singular only at $p=\omega_i$, and $F-F_i$ is analytic at $\omega_i$.
\end{proof}

\begin{Lemma}\label{decaytFi}
  Each function $\mathrm{e}^{-\mu_i p} F_i$ decays like $1/p$ as $p\to\infty$. 
\end{Lemma}
\begin{proof}
  As in Lemma \ref{singFi} 
  , by contour deformation we may assume $p$ lies within the contour of integration defining $F_i$.  Otherwise, by Lemma \ref{singFi} 
  we know the only singularity of of $F_i$ to be $\omega_i$ which lies within the contour and by assumption, $F_i$ is exponentially bounded at infinity.  Using $\eqref{eq:ResidueResurg}$ we see 
  $$p\mathrm{e}^{-\mu_i p} F_i(p)=p\mathrm{e}^{-\mu_i p}F(p)+\frac{1}{2\pi i}\int\limits_{\widetilde{\mathcal{C}}_i} \frac{pF(s)\exp(-\mu_i s)}{s-p}ds$$
  The first term is decaying exponentially. In the second term, we notice that the integrand converges pointwise to the function $-F(s)\mathrm{exp}(-\mu_i s)$ as $p\to \infty$.  Moreover, by assumption $ \mu_is>0$ for large $s$ and $|\mu_i|>\nu$ which imply $$|F(s)\mathrm{exp}(-\mu_i s)|\leq\|F\|_\nu |\mathrm{exp}(-s(\mu_i-\nu))|$$ which is $L^1$.  Using the dominated convergence theorem and the singlevaluedness of the exponential function implies the result.
\end{proof}

\begin{Lemma}\label{GeomLemmaRes}
  On any compact set $K\subset \CC$, the sum in \eqref{eq:decp} converges at least as fast as $\sum_{j\in\ZZ^+,k=1,...,n}\mathrm{e}^{- c_k  j|\lambda_k|(\mu-\nu)}$ (for some $c_k>0$) in the $L^\infty\left(K\right)$ norm to an analytic function on $K$.
\end{Lemma}

\begin{proof}
Let $K\subset \CC$ be compact. Then $K$ only intersects finitely many contours $\mathcal{C}_{1},...,\mathcal{C}_{M}$ (defined in \S\ref{PT}(iv)) and contains (at most) the singular points $\omega_{1},...,\omega_{M}$.  Consider the right side of \eqref{eq:decp}:
$\Phi(p):=F(p)-\sum_i F_i(p)$ for $p\in K$. By Lemma \ref{singFi}, we see $F(p)-\sum_{i=1}^{M} F_i(p)$ is analytic at each $\omega_i,i=1,...,M$ and is therefore analytic throughout $K$.  
Therefore, to prove Lemma\,\ref{GeomLemmaRes} it suffices to do so for the series $\sum_{i>M} F_i(p)$, 

Consider now $i>M$. The contours $\mathcal{C}_i$ are disjoint from $K$ and $F_i$ defined by \eqref{eq:F_i} are analytic for $p\in K$, and are estimated as follows.  
\begin{Remark}\label{RemCpt}
   Note that the following lines leading to estimate \eqref{NonIntEst} do not rely on the compactness of $K$ only that it is a closed set disjoint from $\mathcal{C}_i$.  We mention this because it will be used in the proof of Lemma \ref{FiL1nu} for contours (rays) disjoint from $\mathcal{C}_i$.
\end{Remark} 

Let $\omega_i=j\lambda_k$. First, since $p\in K$ then $p\not\in S_i$ and we can deform the path of integration $\mathcal{C}_i$ in \eqref{eq:F_i} to a path going around the cut $\omega_i+e^{i\theta_k}\RR_+$ (where $\theta_k$ are defined in \S\ref{PT}(iv)). Then, using \eqref{nubounded} and the assumption {\em (A3)},

\begin{multline}\label{NonIntEst}
     \left|F_i(p)\right |\leq\frac{\left|{\rm{e}}^{\mu_i p}\right|}{2\pi d_i}\, 2\, \int_{\omega_i+e^{i\theta_k}\RR_+}\,  \left|F(s)\right | \left|{\rm{e}}^{-\mu_i s}\right|  \left|ds\right| 
     \\
     \leq \frac{\left|{\rm{e}}^{\mu_i p}\right|}{\pi d_i}\, \int_{\omega_i+e^{i\theta_k}\RR_+} \
     \left|F(s)\right| e^{-\nu|s|} \ \left(\sup  \left|e^{-\mu_is}\right|e^{\nu|s|}\right) \left|ds\right|   \le\frac{\left|{\rm{e}}^{\mu_i p}\right|}{\pi d_i}\, \|F\|_{1,\nu}\, e^{ -(\mu-\nu){c}_k j|\lambda_k|}
 \end{multline}

where $d_i=\text{dist}\left(K,\mathcal{C}_i\right)>d_0>0$ and the {\em sup} in \eqref{NonIntEst} is estimated in a straightforward way: for $s=\omega_i+te^{i\theta_k}$ ($t>0$) we have
$\left|e^{-\mu_is}\right|e^{\nu|s|}\le e^{ -(\mu-\nu){c}_k j|\lambda_k|}$ where  ${c}_k=\tfrac{\mu \cos(\theta_k-\arg\lambda_k)-\nu}{\mu-\nu}$ which is positive for $\delta_{2,k}$ small enough. 
Since $\left |\mu_i\right|=\mu$ and $K$ is compact then
$\sup_{i\in\NN}\sup_{p\in K} \left|{\rm{e}}^{\mu_i p}\right|<\infty$.

The estimate in Lemma \ref{GeomLemmaRes} follows. Therefore the right side of \eqref{eq:decp} is analytic on $K$. 

\end{proof}           

\begin{Lemma}\label{FiL1nu}
Assume $F$ satisfies \eqref{nubounded}, $F_i$ are defined by \eqref{eq:F_i}, and use the notations of \S\ref{PT}.
Then there exists $\nu'>\mu$ such that for each singular point $\omega_i$, the corresponding $F_i$ belongs to $L^1_{\nu',\varphi}$ for all $\varphi\in [0,2\pi]$, and so does $\sum_i F_i$.  Moreover, both will have finite $\|\cdot\|_{1,\nu'}$ norm defined by the $\sup$ in \eqref{nubounded}.
\end{Lemma}

\begin{proof}

Fix a singularity $\omega_i\equiv\omega_{jk}=j\lambda_k$ and estimate $F_i(p)$ for $p\in e^{i\phi}\RR_+$.
Suppose, for simplicity, that $\arg\lambda_k=0$ (the proof for other arguments is obtained by a rotation), and that 
$\theta_k>0$ (the proof for $\theta_k<0$ is obtained by symmetry). We distinguish two cases depending on $\phi$: (i) in the first case the Hankel contour and Laplace contour are separated by a positive distance; (ii) in the second, they are not. Fix $\delta>0$ such that ${\rm{e}}^{i(2\pi-\delta)}\RR^+\cap S_i=\emptyset$ and ${\rm{e}}^{i(\theta_k+\delta)}\RR^+\cap S_i=\emptyset$.

(i) For any $\varphi=\arg p\in \left[\theta_k+\delta,2\pi-\delta\right]$ we have $d_{\varphi,\delta}=$ dist$(p, S_i)\ge d_{k,\delta}>0$ , where $d_{k,\delta}$ is the distance from the boundary of the sector $\left\{p:\arg p\in \left[\theta_k+\delta,2\pi-\delta\right]\right\}$ to $S_i$. 
 We now estimate the $L^1_{\nu',\varphi}$ of $F_i$ by splitting the ray into two parts; 
 ${\rm{e}}^{i\varphi}\RR^+={\rm{e}}^{i\varphi}[0,R]\cup{\rm{e}}^{i\varphi}[R,\infty)$ for some $R>0$.  By Lemma \ref{singFi} $\omega_i$ is the only singular point of $F_i$. 
 Furthermore, by analyticity the function $F_i$ is bounded on ${\rm{e}}^{i\varphi}[0,R]$ and hence this part of the norm estimate is finite.  Finally, Lemma \ref{decaytFi} shows that for large $|p|$ the function ${\rm{e}}^{-\nu'|p| }|F_i(p)|$ decays exponentially because $\nu'>\mu$. Therefore integration along ${\rm{e}}^{i\varphi}[R,\infty)$ is finite.  Hence $F_i\in L^1_{\nu',\varphi}$ for all $\varphi\in\left[\theta_k+\delta,2\pi-\delta\right]$.

 (ii) If $\varphi=\arg p\in \left(2\pi-\delta, \theta_k+\delta\right),$  then the two contours may intersect and we can always arrange that they intersect twice by deforming the Hankel contour.  We note that by assumption all singularities $\omega_i$ are locally integrable and so there need not be any distinction between $\varphi\neq 0$ and $\varphi=0$.  We break the Laplace contour into three pieces ${\rm{e}}^{i\varphi}\RR^+ =\gamma_1\cup\gamma_2\cup\gamma_3$. The segments $\gamma_i$ are defined as:
 \begin{equation}\label{GammaDecThreePiec}
     \gamma_1={\rm{e}}^{i\varphi}[0,R_1], \; \gamma_2={\rm{e}}^{i\varphi}(R_1,R_2), \; \gamma_3={\rm{e}}^{i\varphi}[R_2,\infty)
 \end{equation}

where $\gamma_1,\gamma_3$ are a positive distance from $\mathcal{C}_i$, $d_{12}=$dist$\left(\left\{R_1{\rm{e}}^{i\varphi},R_2{\rm{e}}^{i\varphi}\right\},\mathcal{C}_i\right)>0$.  Using characteristic functions we write $F_i=\left(\mathbbm{1}_{\gamma_1}+\mathbbm{1}_{\gamma_2}+ \mathbbm{1}_{\gamma_3}\right)F_i$.  Again estimating the weighted $L^1_{\nu',\varphi}$ norm we see that integration of ${\rm{e}}^{-\nu'|p| }|F_i(p)|$ along both $\gamma_1$ and $\gamma_3$ will be finite by the same argument as in (i) and hence $\left(\mathbbm{1}_{\gamma_1}+ \mathbbm{1}_{\gamma_3}\right)F_i \in L^1_{\nu',\varphi}$.  For any $p\in \gamma_2$ we deform the Hankel contour past $p$ collecting the residues as in \eqref{eq:ResidueResurg} and we arrange the final contour $\widetilde{\mathcal{C}}_i$ so that it satisfies dist$(\gamma_2,\widetilde{\mathcal{C}}_i)\geq d_{12}$; we multiply \eqref{eq:ResidueResurg} by $\mathbbm{1}_{\gamma_2}$ and obtain:
 
  \begin{equation}
        \mathbbm{1}_{\gamma_2}(p)F_i(p)=\mathbbm{1}_{\gamma_2}(p)\left(F(p)+\frac{\exp(\mu_i p)}{2\pi i}\int\limits_{\widetilde{\mathcal{C}}_i} \frac{F(s)\exp(-\mu_i s)}{s-p}ds\right)
  \end{equation}

 Since now $p$ is strictly inside  $\widetilde{\mathcal{C}}_i$, we have $|s-p|\geq d_{12}>0$. The assumptions on $F$ imply that $F_i\mathbbm{1}_{\gamma_2}\in L^1_{\nu',\varphi}$.  Therefore we see $F_i \in L^1_{\nu',\varphi}$ for all $\varphi\in (2\pi-\delta,\theta_k-\delta)$ and by (i) this holds for all $\varphi\in[0,2\pi]$.

 We now consider integrability of the sum.  In the case of (i) $\varphi\in \left[\theta_k+\delta,2\pi-\delta\right]$ we estimate $F_i$ along ${\rm{e}}^{i\varphi}\RR^+$ using \eqref{NonIntEst} (see also Remark \ref{RemCpt}).  This yields:

 \begin{equation}\label{SumEstFi}
     \left|F_i(p)\right|{\rm{e}}^{-\nu'|p|}\leq {\rm{e}}^{-\nu'|p|}\frac{\left|{\rm{e}}^{\mu_i p}\right|}{\pi d_{i,\varphi}}\, \|F\|_{1,\nu}\, e^{ -(\mu-\nu){c}_k j|\lambda_k|}
 \end{equation}
 where $d_{i,\varphi}=$dist$\left({\rm{e}}^{i\varphi}\RR^+,\mathcal{C}_i\right)>0$ for all $\varphi\in \left[\theta_k+\delta,2\pi-\delta\right]$.  Recalling that $F_i\equiv F_{jk}$ and that $\nu'>\mu$ we integrate  \eqref{SumEstFi}. Summing over $j\in \NN$ we see $\sum_{j=1}^\infty \|F_{jk}\|_{1,\nu',\varphi}<\infty$.  Since $L^1_{\nu',\varphi}$ is a Banach space this implies $\sum_{j=1}^\infty F_{jk}\in L^1_{\nu',\varphi}$ and the same is true for  $\sum_i F_i=\sum_{k=1}^n\sum_{j=1}^\infty F_{jk}$.
 Next, we now consider $\varphi\in \left(2\pi-\delta, \theta_k+\delta\right)$ and show the sum is integrable along these directions.  We proceed as in (ii) and arrange that each Hankel contour involved in the sum intersects the Laplace contour exactly twice.  We then decompose each Laplace contour as we did in \eqref{GammaDecThreePiec}.  More precisely, for each $\mathcal{C}_i\equiv \mathcal{C}_{jk}$:\,${\rm{e}}^{i\varphi}\RR^+ =\gamma_{j_1}\cup\gamma_{j_2}\cup\gamma_{j_3}$ where $\gamma_{i_j}$ are defined as:
 \begin{equation}\label{MultGammaDecThreePiec}
     \gamma_{j_1}={\rm{e}}^{i\varphi}[0,R_{j_1}], \; \gamma_{j_2}={\rm{e}}^{i\varphi}(R_{j_1},R_{j_2}), \; \gamma_{j_3}={\rm{e}}^{i\varphi}[R_{j_2},\infty)
 \end{equation}
where we place the contours $\gamma_{j_1},\gamma_{j_3}$ at a small enough positive distance from $\mathcal{C}_i$; $d=$dist$\left(\left\{\gamma_{j_1},\gamma_{j_3}\right\},\mathcal{C}_{jk}\right)>0$. This is possible because of equal separation of singularities along each singular array.
Estimating the sum of the norms we have 

\begin{equation}\label{NormofSumEstSingDir}
    \sum_{j=1}^\infty \int_{{\rm{e}}^{i\varphi}\RR^+} e^{-\nu'|p| }\,\left| F_{jk}(p)\right|\, |dp|= \sum_{j=1}^\infty\left( \int_{0}^{R_{j_1}}+ \int_{R_{j_1}}^{R_{j_2}} + \int_{R_{j_2}}^\infty \right) {\rm{e}}^{-\nu' p }\,\left| F_{jk}(p{\rm{e}}^{i\varphi})\right|\, dp
\end{equation}
Using the estimate \eqref{NonIntEst} we get an upper bound on the first and third inetgrals

\begin{multline}
        \left(\int_{0}^{R_{j_1}}+ \int_{R_{j_2}}^\infty \right) {\rm{e}}^{-\nu' p }\,\left| F_{jk}(p{\rm{e}}^{i\varphi})\right|\, dp \\
        \leq 
        \frac{{\rm{e}}^{-(\mu-\nu)c_k j |\lambda_k|}}{\pi d}\|F\|_{1,\nu}\left(\int_{0}^{R_{j_1}}+ \int_{R_{j_2}}^\infty \right){\rm{e}}^{-\nu' p }\, {\rm{e}}^{|\mu_i|p\cos(\varphi+\arg \mu_i)} dp\\
        \leq \frac{{\rm{e}}^{-(\mu-\nu)c_k j |\lambda_k|}} {\pi d}\|F\|_{1,\nu}\int_0^\infty {\rm{e}}^{-p(\nu' -\mu \cos(\varphi+\arg \mu_i)} dp=\frac{{\rm{e}}^{-(\mu-\nu)c_k j |\lambda_k|}} {\pi d(\nu' -\mu \cos(\varphi+\arg \mu_i)}\|F\|_{1,\nu}
\end{multline}
By assumption $\mu>\nu$ , $\nu'>\mu$ and ${c}_k=\tfrac{\mu \cos(\theta_k-\arg\lambda_k)-\nu}{\mu-\nu}>0$  for $\delta_{2,k}$ small enough; see \eqref{assumethetak}. Hence, the first and the third sums in \eqref{NormofSumEstSingDir}  converge.

If $p\in \gamma_{j_2}$ we once more deform the contour $\mathcal{C}_{jk}$ past $p$ collecting the residues and we arrange the final contour $\widetilde{\mathcal{C}}_{jk}$ so that it satisfies dist$(\gamma_{j_2},\widetilde{\mathcal{C}}_{jk})=d>0$.

We use \eqref{eq:ResidueResurg} with $|p-s|\ge d>0$ for $s\in \widetilde{\mathcal{C}}_{jk}$ and $p\in \gamma_{j_2}$.   

\begin{multline}\label{SumEstSing}
     \sum_{j=1}^\infty  \int_{R_{j_1}}^{R_{j_2}}    {\rm{e}}^{-\nu' p }\,\left| F_{jk}(p{\rm{e}}^{i\varphi})\right|\, dp \\
   \leq \sum_{j=1}^\infty  \int_{R_{j_1}}^{R_{j_2}}    {\rm{e}}^{-\nu' p }\,\left( \left| F(p{\rm{e}}^{i\varphi})\right|+\frac{|\exp(\mu_{jk} p{\rm{e}}^{i \varphi })|}{2\pi }\,\int\limits_{\widetilde{\mathcal{C}}_{jk}} \left|\frac{F(s)\exp(-\mu_{jk} s)}{s-p{\rm{e}}^{i \varphi }}\right||ds|\right)\, dp 
\end{multline}

Using integrability $F\in L_\nu^1$ and $\nu'>\nu$ we have

\begin{equation}
    \sum_{j=1}^\infty  \int_{R_{j_1}}^{R_{j_2}}    {\rm{e}}^{-\nu' p }\, \left| F(p{\rm{e}}^{i\varphi})\right|dp\leq \int_{{\rm{e}}^{i \varphi }\RR^+}{\rm{e}}^{-\nu' p }\, \left| F(p{\rm{e}}^{i\varphi})\right|dp\le \|F\|_{1,\nu}<\infty
\end{equation}

Once again we estimate using \eqref{NonIntEst} and see the last term in \eqref{SumEstSing} is finite:

\begin{multline}
     \sum_{j=1}^\infty  \int_{R_{j_1}}^{R_{j_2}}    {\rm{e}}^{-\nu' p } \frac{|\exp(\mu_{jk} p {\rm{e}}^{i \varphi })|}{2\pi }\int\limits_{\widetilde{\mathcal{C}}_{jk}} \left|\frac{F(s)\exp(-\mu_{jk} s)}{s-p{\rm{e}}^{i \varphi }}\right||ds| \, dp \\
     \leq \frac{\|F\|_{1,\nu}}{\pi d}\sum_{j=1}^\infty {\rm{e}}^{-(\mu-\nu)c_k j |\lambda_k|}\int_{R_{j_1}}^{R_{j_2}}{\rm{e}}^{-p(\nu'-\mu\cos(\arg\mu_{jk}+\varphi))  }dp\\
     \leq\frac{\|F\|_{1,\nu}}{\pi d(\nu'-\mu\cos(\arg\mu_{jk}+\varphi))}\sum_{j=1}^\infty {\rm{e}}^{-(\mu-\nu)c_k j |\lambda_k|}<\infty
\end{multline}
Therefore we see $\sum_{j=1}^\infty \|F_{jk}\|_{1,\nu',\varphi}<\infty$ if $\nu'>\mu$ which implies $\sum_iF_i\in L_{\nu',\varphi}^1$.  Hence we conclude that both $F_i$ and $ \sum_iF_i\in L_{\nu',\varphi}^1$ for every $\varphi\in[0,2\pi]$.  

Finally we show $F_i$ and $ \sum_iF_i$ have finite  $\|\cdot\|_{1,\nu'}$ norm.  We show this first for $F_i$.    Recalling the notation $\omega_i\equiv \omega_{jk}=j\lambda_k$, given any $\gamma\in \mathscr{R}$ choose a representative so that dist$(\gamma,\{j\lambda_k\}_{j\in\NN})=d>0$. Again, we use \eqref{NonIntEst} applied to each point $p\in\gamma$ and to each resurgent element $F_{jk}$ and obtain the uniform bound (see also \eqref{nubounded})

\begin{multline}\label{GamIntEst}
     \int_{\gamma} {\rm{e}}^{-\nu'|p|}|F_{jk}(p)||dp|\leq \frac{1}{\pi d}\|F\|_{1,\nu}\, e^{ -(\mu-\nu){c}_k j|\lambda_k|}\int_\gamma {\rm{e}}^{-\nu'|p|}\left|{\rm{e}}^{\mu_{i} p}\right||dp|\\
     \leq \frac{K}{\pi d}\|F\|_{1,\nu}\, e^{ -(\mu-\nu){c}_k j|\lambda_k|}, \quad \text{where} \, K=\sup _{\gamma\in\mathscr{R}}\int_\gamma {\rm{e}}^{-(\nu'-\mu)|p|}|dp|<\infty
\end{multline}
by \eqref{defR}.  Therefore we have $F_{jk}\in L_{\nu'}^1$, $\sum_j \|F_{jk}\|_{1,\nu'}<\infty$ and it follows that $\sum_jF_{jk}\in L^1_{\nu'}$.  The same holds for $\sum_{k=1}^n \sum_jF_{jk}\in L^1_{\nu'}$.

\end{proof}

\begin{Lemma}\label{GisLaplace}

Let $\nu'$ be as in Lemma\,\ref{FiL1nu}.

  The function
  \begin{equation}
    \label{eq:eqG}
    G(p)=F(p)-\sum_{i}F_i
  \end{equation}
is entire and $G(p)\mathrm{e}^{-\mu'|p|}\in L^{\infty}(\mathrm{e}^{i\varphi}\RR^+)$ for any $\mu'\geq\nu'$ and any $\varphi\in[0,2\pi]$.
\end{Lemma}
\begin{proof}
Note that $G\in L^1_{\nu',\varphi}$ by Lemma\,\ref{FiL1nu} and assumption \eqref{nubounded}. 

  Analyticity follows from the monodromy theorem, since $G$ has analytic continuation along any ray in $\CC$. To see the bound we consider the antiderivative of $G$.  Let $H(p)=\int_0^pG(s)ds$ and estimate 
  $$ |H(p)|\leq \mathrm{e}^{\nu'|p|} \int_0^p \mathrm{e}^{-\nu'|s|}\left|G(s)\right|\,|ds| \leq\mathrm{e}^{\nu'|p|}\|G\|_{1,\nu'}$$
  hence $H$ is of exponential order one.  Let $R>0$, using Cauchy's theorem and $G(p)=H'(p)$ we have
  $$|G(p)|\leq \frac{1}{2\pi }\oint_{|s|=R}\frac{|H(s)|}{|s-p|^2}|ds|\leq \frac{1}{2\pi }\int_{0}^{2\pi}\frac{\mathrm{e}^{\nu'|p|}\|G\|_{1,\nu'}}{(R-|p|)^2}Rdt =\frac{R\mathrm{e}^{\nu'|p|}\|G\|_{1,\nu'}}{(R-|p|)^2}$$
  for  any $|p|<R$.  Letting $R\to \infty$ implies the uniform exponential bound.

\end{proof}
\begin{Lemma}
  $g=\mathcal{L}G$ has a convergent asymptotic series at infinity, and is equal to the sum of the series. 
\end{Lemma}
\begin{proof}
  Expanding $G$ into a power series about the origin $G(p)=\sum_{n=0}^\infty a_n p^n$ with $a_n=\frac{G^{(n)}(0)}{n!}$ which is also an asymptotic series.  Watson's lemma shows 
  $$g=\mathcal{L}G\sim \sum_{n=0}^\infty \frac{n!a_n}{x^{n+1}}=\sum_{n=0}^\infty \frac{G^{(n)}(0)}{x^{n+1}}$$
  as $x\to \infty$ along any direction in $\CC$.\\
  Moreover, the uniform bound we have for $G$ and Cauchy estimates provide $|G^{(n)}(0)|\lesssim \mu'^n$ where $\mu'$ is as in \eqref{eq:eqG}.  Therefore, if $|x|$ is large enough the expansion for $g$ at infinity will converge.  The function $h(z)=g(1/z)$ is bounded at zero and single-valued, as is seen by deformation of contour (since $G$ is exponentially bounded and entire). Thus $h$ is analytic at zero, and therefore the sum of its asymptotic (=Taylor) series at zero.
\end{proof}

\begin{Lemma}\label{finlem}
  The change of variable $\tilde{x}=x-\mu_i$ leads to
$\lap[F_i](x)=\lap[\tilde{F}_i](\tilde{x})$ where $\tilde{F}_i$ decays like $1/p$ as $p\to\infty$. 
\end{Lemma}
Combining these lemmas, Theorem\,\ref{TT1} follows. 

\end{proof}
 
\section{Appendix}

\subsection{Normalized Airy and Bessel functions}\label{sec5}
The modified Bessel equation is 
\begin{equation}
  \label{eq:eqb}
  x^2y'' +xy'-(\nu^2+x^2)y=0
\end{equation}
The transformation 
$y=\mathrm{e}^{-x}x^{1/2}h(2x)$, $u=2x$ brings \eqref{eq:eqb} to the normalized form  
\begin{equation}
  \label{eq:eqnormB}
  h''-\left(1-\frac2u\right)h'-\left(\frac1u-\frac1{4u^2}+\frac{\nu^2}{u^2}\right) h=0
\end{equation}
This normalized form is suitable for Borel summation since it admits a formal power series solution in powers of $u^{-1}$ starting with $u^{-1}$; it is further normalized to ensure that the Borel plane singularity is placed at $p=-1$. One way to obtain the transformation is to rely on the classical asymptotic behavior of Bessel functions and seek a transformation that formally leads to a solution as above. 

The Airy equation 
\begin{equation}
  \label{eq:Ai1}
 f''-xf=0
\end{equation}
can be brought to the Bessel equation with $\nu=1/3$. The normalizing transformation can be obtained from \eqref{eq:eqnormB}, or directly, based on the asymptotic behavior at $\infty$ which suggests the change of variables
$$f(x)=x^{5/4}\mathrm{e}^{-\frac23 x^{3/2}}h(4x^{3/2}/3);\ \ \ \ \ \ x=(3u/4)^{2/3}$$
bringing the equation to
\begin{equation}
  \label{eq:eqa1}
  h''-\left(1-\frac2{u}\right)h'-\left(\frac1{u}-\frac{5}{36u^2}\right)h=0
\end{equation}
which is indeed \eqref{eq:eqnormB} for $\nu=1/3$. From this point, without notable algebraic complications we analyze \eqref{eq:eqnormB}.

The inverse Laplace transform of \eqref{eq:eqnormB} results in an integral equation which we differentiate twice to obtain:
\begin{equation}
  \label{eq:borel}
  p (p+1) H''(p)+(2 p+1) H'(p)+\left(\frac{1}{4}-\nu ^2\right) H(p) =0
\end{equation}
which is a hypergeometric equation \cite{nist}15.10. Its solution which is {\em analytic} at $p=0$ is (a constant multiple of)
\begin{equation}
  \label{eq:solna}
 \ _2F_1\left(\tfrac12+\nu,\tfrac12-\nu;1;-p\right)=P_{\nu-\frac12}(1+2p)
\end{equation}
where $\,_2F_1$ is the usual hypergeometric function \cite{nist}(14.3.1) and $P_\mu$ is the Legendre $P$ function. On the first Riemann sheet, the solution has two regular singularities, $p=-1$ and $p=\infty$. The behavior at zero is \cite{nist}(15.2.1)
$$ P_{\nu-\frac12}(1+2p)=1+\left(\nu ^2-\frac{1}{4}\right) p+\frac{1}{64} \left(16 \nu ^4-40 \nu
   ^2+9\right) p^2+\cdots$$

At $p=-1$ (with the phase of the log being the usual one for $p+1>0$ we have
\begin{equation}
  \label{asympism1}
  P_{\nu-\frac12}(1+2p)=c_1\,\log (p+1)+c_2+c_3 (p+1)\log(p+1)+c_4(p+1)+\ldots
  \end{equation}
 where $c_j=c_j(\nu)$ are nonzero constants.

 Using properties of the hypergeometric function \cite{nist}(15.2.3) $_2F_1(a,b;c;z)$ and its relation \eqref{eq:solna} to the Legendre $P$ function, we recover the branch jump across the cut $(-\infty,-1]$ 
 
 \begin{equation}
     \label{eq:Pbranch}
 P_{\nu-\frac{1}{2}}^+(1+2p)-P_{\nu-\frac{1}{2}}^-(1+2p)=\frac{2\pi i}{\Gamma\big(\frac{1}{2}+\nu\big)\Gamma\big(\frac{1}{2}-\nu\big)}P_{\nu-\frac{1}{2}}(-1-2p)=2i\cos (\pi \nu)P_{\nu-\frac{1}{2}}(-1-2p)
 \end{equation}
Where we used the reflection formula for the Gamma function to obtain the last equality. 
We differentiate \eqref{eq:Pbranch} to obtain the branch jump for the derivative which will be used in \S\ref{Airy}.  
   
\begin{equation}
    \label{eq:derPbranch}
   \Delta  P_{\nu-\frac{1}{2}}'(1+2p)= -2i\cos (\pi \nu)P_{\nu-\frac{1}{2}}'(-1-2p)
\end{equation}
We note that such a simple relation for the branch jump stems from the fact that the Airy function satisfies a linear second order ODE and does not hold in general.

  \subsection{The derivatives of the polylogarithm}\label{form} Here we show that the coefficients  $s(k,j)$ which appear in \eqref{Llicoeffs} are the Stirling numbers of the first kind.
  
  For $k=0$ we have $s(0,0)=1$. It is easy to check that Li$_s'(z)=z^{-1}$Li$_{s-1}(z)$ confirming that $s(1,0)=0$ and $s(1,1)=1$. For higher $k$ formula \eqref{derLi} is then checked by a simple induction, which leads to the recurrence relations
 $$s(k+1,0)=-ks(k,0),\ \ s(k+1,k+1)=s(k,k),\ \ s(k+1,j)=-ks(k,j)+s(k,j-1)$$
 which are the recurrence relations satisfied by the Stirling numbers of the first kind, see \cite{nist} Sec.26.8.
 
 We note that the polylog is another example when the factorial series converges geometrically in a cut plane.\footnote{We note that ${\rm Li}_s(z)=z\Phi(z,s,1)$.}

  \subsection{More about Ei}\label{Eidetails}
  
 We note that $E_1$ can be written in a form that allows for analytic continuation through the cut on $\RR_-$: by elementary changes of variables we have, for $z>0$,
  $$E_1(z):=\int_z^\infty\frac{e^{-t}}{t}\, dt= e^{-z}\int_0^\infty \frac{e^{-pz}}{1+p}\, dp$$
 
 It is interesting to note the relation \cite{nist} 6.2.4
 $$E_1(z)={\rm Ein}(z)-\ln z-\gamma,\ \ \ \ \ \ \text{where }\ \ \ {\rm Ein}(z)=\int_0^z\frac{1-e^{-t}}{t}\, dt$$
 Since Ein is an entire function, and ln is defined with the usual branch for $z>0$, we see that upon analytic continuation across $\RR_-$, the function $E_1$ gains a $2\pi i$; thus $\RR_-$ is a Stokes line.
 
 For us it is convenient to place this Stokes line along $\RR_+$, so we work with $-E_1(-z)$:
 $$-E_1(-z)= -e^{z}\int_0^\infty \frac{e^{pz}}{1+p}\, dp,\ \ (z<0)$$
 which analytically continued to the first quadrant yields our Ei$^+$ defined in \eqref{LiEi}.
 
 Note the structure of the branch point at $0$:
 
 $${\rm Ei}^+(z)=-{\rm Ein}(-z)+\ln (-z)+\gamma$$
 where $\ln (-z)$ has the usual brach for $z<0$ and then it is analytically continued on the Riemann surface of the log (and ${\rm Ein}(-z)$ is entire).
 
 For $e^{-x}$Ei$^+(x)$, $\RR^+$ is a Stokes ray and the two sides of $i\RR^-$ are antistokes lines. The behavior of $e^{-x}$Ei$^+(x)$ is oscillatory when $\arg(x)=-\pi/2$ and it is given by an asymptotic series when $\arg(x)=3\pi/2$.

 \subsection{Proof of Lemma\,\ref{RemLer}}\label{ProofL11}

 The function $\Phi(z,1,x)$ is defined as (see \cite{nist} 25.14.1)
 $$\Phi(z,1,x)=\sum_{n=0}^\infty\frac{z^n}{x+n},\ \ \ \ \ \text{for } x\in\CC\setminus\ZZ_- ,\  |z|<1$$
 and for other values of $z$, it is defined by analytic continuation.
 
 Clearly $\Phi(z,1,x)$ is a meromorphic function of $x$ (for $|z|<1$).

 For $\Re x>0$ and $|z|<1$, $\Phi$ has the integral representation 
 \begin{equation}
 \label{IntPhi}
 \Phi(z,1,x)=\int_0^\infty\frac{e^{-xp}}{1-ze^{-p}}dp
 \end{equation}
  (see \cite{nist} 25.14.5). Note that the right hand side of \eqref{IntPhi} is analytic in $z$ in $\CC\setminus [1,\infty)$, hence  \eqref{IntPhi} holds in this domain.
  
  Then, for $\Re x>0$ and $|z|<1$, by dominated convergence and using \eqref{eq:identexp}, we get
   \begin{multline}
     \label{eq:exp-expansion}
 \Phi\left( \frac{z}{z-1},1,x\right)=(1-z)\int_0^\infty \frac{e^{-xp}}{1-z(1-e^{-p})}dp=(1-z)\int_0^\infty e^{-xp}\sum_{k\ge 0}z^k(1-e^{-p})^kdp\\
 =(1-z)\sum_{k\ge 0}z^k\mathcal{L}(1-e^{-p})^k=(1-z)\sum_{k\ge 0}z^k\frac{k!}{(x)_{k+1}}:=(1-z)\mathcal{E}(x,z)
   \end{multline}
  We now show that the series $\mathcal{E}(x,z)$ converges uniformly for $x$ in any compact set disjoint from $\ZZ_-$. Once this is proved it follows that the equality 
  \begin{equation}
     \label{expPhiseries}
 \Phi\left( \frac{z}{z-1},1,x\right)=(1-z)\sum_{k\ge 0}z^k\frac{k!}{(x)_{k+1}}
   \end{equation}
    holds for all $x\in\CC\setminus\ZZ_- ,\  |z|<1$.
  
  \
  
We fix $z$ with $|z|<1$.  Let $K$ be a compact set disjoint from $\ZZ_-$.  Let $\lambda$ be such that $|z|\le\lambda<1$. Let $M$ be as assumed in \eqref{valM}.

We prove that $S_n=\sum_{k= M}^n z^k\frac{k!}{(x)_{k+1}}$ converges uniformly in the following steps outlined below.

{\em (i)}  Clearly $S_n$ satisfies the recurrence 
\begin{equation}
\label{recf}
f_n-f_{n-1}=z^n\frac{n!}{(x)_{n+1}}
\end{equation}
 {\em (ii)} We show that the recurrence \eqref{recf} has a unique solution of the form 
 \begin{equation}
 \label{fnvn}
 f_n=z^n\frac{n!}{(x)_{n+1}}\left(\frac{z}{z-1}+\frac{u_n}{1-z}\right)
 \end{equation}
  where $u_n\to 0$. {\em (iii)} We then show that $f_n\to 0$. 
{\em (iv)} Since the general solution of the recurrence \eqref{recf} is $f_n+C$, it follows that  $\sum_{k= M}^n z^k\frac{k!}{(x)_{k+1}}=f_n+C$ for some $C$, hence the series $\sum_{k= M}^\infty z^k\frac{k!}{(x)_{k+1}}$ converges. {\em (v)} Finally, we estimate the remainder of this series.

\

{\em Proof of (ii)-(iv).}

{\em (ii)} Assume $\{f_n\}_n$ satisfies \eqref{recf}. Let $u_n$ be given by \eqref{fnvn}. Then $u_n$ satisfies
 \begin{equation}
 \label{recvn}
 u_{n-1}=\frac{xz}{(x+n)}+\frac{nz}{x+n}u_n:=a_{n-1}+d_{n-1}\,u_nx
 \end{equation}
 which,  for $n\ge M+1$ can be written as a functional equation: denoting ${ \bf u}=\langle u_{M},u_{M+1},\ldots\rangle $ the recurrence \eqref{recvn} is 
 \begin{equation}
 \label{equ}
 { \bf u}={ \bf a}+DS{ \bf u}
  \end{equation}
   where $S$ is the left shift and $D$ is a diagonal operator. We consider this equation in a weighted $\ell^\infty$: let $\mathcal{B}=\{{ \bf u}\, |\, \|{ \bf u}\|=\sup_{n\ge M-1} n|u_n|<\infty\}$.

 We show that $\|DS\|\le\lambda <1$ and then apply the contractive mapping principle.
   Indeed, for each ${ \bf u}\in\mathcal{B}$ 
  \begin{equation}
 \label{estiml}
 |(n-1)(DS{ \bf u})_{n-1}|=(n-1)\frac{n|z|}{|x+n|}|u_n|\le  \|{\bf u}\| \frac{|z|}{\left|\frac xn+1\right|}
 \end{equation}
 
 Note that 
 $\left|\frac xn+1\right|^2=\left|\frac xn\right|^2+2\left|\frac xn\right|\cos\alpha +1$.

 We show that $\left|\frac xn+1\right|^2>|z|^2/\lambda^2$ hence, by \eqref{estiml}, $\|DS\|\le\lambda$, and thus $DS$ is a contraction. We see that if $\Re x\ge 0$ then $\left|\frac xn+1\right|^2\ge 1>|z|^2/\lambda^2 $ therefore it suffices to consider the case $\Re x=|x|\cos \alpha<0$.  Denote $\zeta=i|z|/\lambda$ and $\xi=\Re x/|x|$. If $|\zeta-\xi|<1$ we can see that $\left|\frac xn+1\right|^2\ge 1$ for any $n\ge 1$. However, if $|\zeta-\xi|\ge 1$ we impose the condition $n> M$ with $M$ as in \eqref{valM}.\footnote{ Indeed, assuming $n>M$, $|\zeta-\xi|^2=\cos^2\alpha+|z|^2/\lambda^2\ge 1$ and  taking $M$ to be as in \eqref{valM}  we again have the desired bound. }

 Therefore $\|DS\|\le\lambda<1$.
  
 The value of $M$ is as follows. $M=1$ if $\Re x<0$ and $x\in B_1(z):=\{\zeta:|\zeta-z|<1\}$. Otherwise $M$ is given by  \eqref{valM}.
 
 Also  
 \begin{equation}
 \label{estima}
 \|\mathbf{a}\|=\sup_{n\ge M}(n-1)|a_{n-1}|= |xz|\sup _{n\ge M} \frac{n-1}{|x+n|} \le |x| \lambda
 \end{equation}
 therefore $\mathbf{a}\in \mathcal{B}$.
 
 Therefore this contractive linear equation has a unique solution in $\mathcal{B}$:  ${ \bf u}=(I-DS)^{-1}{ \bf a}$. 
 
 This solution satisfies
 $\|\mathbf{{u}} \|\leq \tfrac1{1-\lambda}\|\mathbf a\|$ hence, using \eqref{estima}
 $$|{u_n}|\le \frac 1n \|\mathbf{{u}}\|\le \frac 1n \frac1{1-\lambda} |xz|\, \sup _{n\ge M} \frac{n-1}{|x+n|} $$

{\em (iii)} To show that $f_n\to 0$ note that, by Stirling's formula, we have, for large $n$ and $x\notin \ZZ_-$ 
 \begin{equation}
   \label{eq:factorial-asympt}
  \frac{n!}{(x)_{n+1}}= \Gamma(x)\frac{n!}{\Gamma(x+n+1)}=\Gamma(x)\,n^{-x}\,(1+\epsilon(n,x)),\ \ \ \ \text{with } |\epsilon(n,x)|<\frac{C}{n}\text{ for }x\text{ in a compact set}  
 \end{equation}
 which goes to $0$ uniformly for $x\in K$. Then {\em (iv)} follows and \eqref{expPhiseries}
     holds for all $x\in\CC\setminus\ZZ_- ,\  |z|<1$.
 
 \
 
 {\em Estimate of the remainder}
 
 {\em (v)} To estimate the remainder, since $\sum_{k= M}^n z^k\frac{k!}{(x)_{k+1}}=f_n+C$ then (for $n\ge M$) we have  $\sum_{k= 0}^n z^k\frac{k!}{(x)_{k+1}}=f_n+C_1$ and since $f_n\to 0$ then $C_1=\sum_{k= 0}^\infty z^k\frac{k!}{(x)_{k+1}}$.  
 
 It follows that $\sum_{k= n+1}^\infty z^k\frac{k!}{(x)_{k+1}}=-f_n$. Since the remainder is $\rho_{n+1,0}
 =-(1-z)f_n$, using \eqref{fnvn} we obtain \eqref{rhonestimate}.
 
  Formula \eqref{newrhonestimate} is obtained integrating by parts $n$ times the integral representation \eqref{eq:exp-expansion}: for  $\Re (x)>0$ we have 
 
 $$\int_0^\infty \frac{e^{-xp}}{1-z(1-e^{-p})}dp=\int_0^1\frac{t^{x-1}}{1-z+zt}\, dt=\sum_{j=0}^n\frac{j!z^j}{(x)_{j+1}}+
z^{n+1}\frac{(n+1)!}{(x)_{n+1}}\, \int_0^{\infty}\frac{ e^{-p(x+n+1)}}{(1-z+ze^{-p})^{n+2}}\, dp$$
 for all $n\ge 0$.
 
 For $\Re x\le 0$ we do analytic continuation by rotating the path of integration in \eqref{eq:exp-expansion} as follows. Note that the points $p$ for which $1-z+ze^{-p}=0$ lie on a vertical line with with $\Re p\le -\ln\left(\tfrac 1{|z|}-1\right):=a_z$. Then we consider a path of integration starting at $0$, going along $\RR_+$ up to a point $p_z>a_z$ and we rotate the rest of the half-line by clockwise, or counterclockwise, at most $\pi/2$ in such a way that $\Re(xp)>0$. 
 
To prove \eqref{newrhonestimate2} we estimate the integral in \eqref{newrhonestimate} using the method of steepest descent. We have, denoting $m=n+1$ ,
$$ \int_0^{\infty }\frac{ e^{-p(x+m)}}{(1-z+ze^{-p})^{m+1}}\, dp=\int_0^{\infty } e^{-m(p+\ln(1-z+ze^{-p}))}\frac{e^{-px}}{1-z+ze^{-p}}\, dp:= \int_0^{\infty } e^{-mf(p)}g(p)dp$$

Noting that $f'(p)\ne 0$, then $\Re f(p)$ has no max/min; thus it is increasing and the dominant behavior is obtained at $p=0$ and the dominant behavior of the integral is $ \int_0^{\infty } e^{-m[f(0)+f'(0)p]}g(0)dp=\tfrac1{m(1-z)}$.

Finally, using Stirling's formula we have $\tfrac{m!}{(x)_m}\sim m^{-x}\Gamma(x)$, and thus \eqref{newrhonestimate} implies \eqref{newrhonestimate2}.

 \section{Acknowledgments} We are very grateful for the referees' extremely thorough reading of the paper and their many useful suggestions.
OC was partially supported by the NSF grants DMS-1515755 and DMS-
2206241.

\end{document}